\documentclass[10pt]{IEEEtran}

\usepackage{amsmath,amsthm}
\usepackage{amssymb}
\usepackage{amsfonts}
\usepackage{fancyhdr}
\usepackage{afterpage}
\usepackage{float}
\usepackage{bm,balance}
\usepackage{algorithm,algorithmic}
\usepackage{graphicx,color,epsfig,psfrag}
\usepackage{multicol,multirow}

\usepackage{cite}      

\usepackage{graphicx}  

\usepackage{subfigure} 

\usepackage{amsmath}   
\usepackage{array}
\theoremstyle{definition}

\theoremstyle{remark}

\newcommand{\btheta}{{\mbox{$\bm{\theta}$}}}

\begin{document}
%
\title{Mean Square Error bounds for parameter estimation under model misspecification}

\author{\authorblockN{
\authorblockN{Adri{\`a} Gusi-Amig{\'o}, \IEEEmembership{Student Member,~IEEE}, Pau Closas, \IEEEmembership{Senior Member,~IEEE} \\ and Luc Vandendorpe, \IEEEmembership{Fellow,~IEEE}.} \\ 
\thanks{Copyright (c) 2014 IEEE. Personal use of this material is permitted. However, permission to use this material for any other purposes must be obtained from the IEEE by sending a request to pubs-permissions@ieee.org}
\thanks{This work has been partially supported by FRIA and the FNRS, by the IAP BESTCOM network, the Government of Catalonia under Grant 2014--SGR--1567, and by the European Commission in the Network of
Excellence in Wireless COMmunications NEWCOM$\sharp$ (contract n.
318306)}
\thanks{A. Gusi-Amig\'o was with the Institute of Information and Communication Technologies, Electronics and Applied Mathematics (ICTEAM), Universit{\'e} catholique de Louvain (Belgium). Now he is with the Centre Tecnol\`{o}gic de Telecomunicacions de Catalunya (CTTC), Parc Mediterrani de la Tecnologia, Av. Carl Friedrich Gauss 7, 08860 Castelldefels, Barcelona (Spain). e-mail: adria.gusi@cttc.cat}
\thanks{P. Closas is with the Centre Tecnol\`{o}gic de Telecomunicacions de Catalunya (CTTC), Parc Mediterrani de la Tecnologia, Av. Carl Friedrich Gauss 7, 08860 Castelldefels, Barcelona (Spain). e-mail: pau.closas@cttc.cat}
\thanks{L. Vandendorpe is with the Institute of Information and Communication Technologies, Electronics and Applied Mathematics (ICTEAM), Universit{\'e} catholique de Louvain, Place du Levant, 2, 1348 Louvain-la Neuve (Belgium). e-mail: luc.vandendorpe@uclouvain.be}
}}
%
%
%
\markboth{submitted to the IEEE Trans. on Signal Processing}{Gusi, Closas, and Vandendorpe: Mean Square Error bounds for parameter estimation under model misspecification}

%



\maketitle

\begin{abstract}
In parameter estimation, assumptions about the model are typically
considered which allow us to build optimal estimation methods under
many statistical senses. However, it is usually the case where such
models are inaccurately known or not capturing the complexity of the
observed phenomenon. A natural question arises to whether we can
find fundamental estimation bounds under model mismatches. This
paper derives a general bound on the mean square error (MSE)
following the Ziv-Zakai methodology for the widely used additive
Gaussian model. The general result accounts for erroneous
functionals, hyperparameters, and distributions differing from the
Gaussian. The result is then particularized to gain some insight
into specific problems and some illustrative examples demonstrate
the predictive capabilities of the bound.
\end{abstract}

\IEEEpeerreviewmaketitle

\section{Introduction}
\label{sec:intro}

\PARstart{M}{odels} are parsimonious representations of nature that
try to capture the most relevant features of a process. To some
extent all models are erroneous since natural phenomena are
typically much more complex than the assumptions typically imposed.
For instance, the Gaussian assumption is typically considered in the
statistical signal processing literature advocating for the central
limit theorem and its interesting mathematical properties
\cite{Kim08}. However, there are plenty of situations where
Gaussianity does not hold, for instance due to the presence of
outliers in the measurements. Refer to the enlightening article
\cite{Zoubir12} and the references therein for examples of the
latter. Nevertheless, the Gaussian model is of paramount
significance in estimation theory, allowing in many cases the
development of realizable algorithms. Therefore, the goal is to
understand when this model is \emph{good} enough for parameter
estimation purposes \cite{Box86}.

Robust statistics is the discipline of statistics that studies the
impact of model mismatch in inference procedures and investigates
how to circumvent this limitation. The foundations of modern robust
statistics can be traced back to the 1960's, where the first
definitions and studies of robust estimators appeared, although
robust statistics emerged as a hot research trend in the mid 1980's
with the first books covering the topic \cite{Huber81,Hampel86}.
Important contributions were published since then, both in the
development of metrics for the assessment of robustness of
estimation methods and in the design of robust methods
\cite{Kassam85,Olive98,Maronna06,Dietrich10}. The achievements in
the latter having its seminal work in Huber's M-estimators
\cite{Huber81}, a broad class of estimators generalizing the Maximum
Likelihood (ML) principle which are designed to be robust to model
departures.

In robust statistics, the focus in robustness assessment has been in
quantifying and interpreting the sensitivity of a method to model
uncertainties. From a classical perspective, robustness can be
defined as the insensitivity to small deviations from the
assumptions. Assuming a parametric model, a statistical procedure
should possess the following features to be considered robust: $i$)
have optimal (or nearly optimal) performance under the assumed model;
$ii$) small deviations from the model assumptions should not have a
large impact on the method's performance (this is referred to as
Qualitative robustness of a method); and $iii$) large deviations
from the model assumptions should not cause a catastrophe (which is
the standard definition of Quantitative robustness). The three
features are typically characterized respectively by: the relative
efficiency, which is the ratio of variances of the optimal method
and the method under study in nominal model conditions; the
stability, measured by the influence function (IF) as the bias
impact of infinitesimal contaminations of the model; and the
breakdown point (BP), defined as the maximal fraction of outliers
the method can handle without collapsing. All these definitions and
metrics are well-established and used in the context of estimation
theory. They have the characteristic of being specific to each
particular estimator, but not general for the estimation problem
under misspecified models.

When deriving estimators for a particular problem, one is typically
interested in the fundamental estimation bounds that can be
achieved. The goal is to evaluate the ability of the estimator to
attain the bound (and thus being efficient). Classical (i.e.,
non-robust) estimation theory provides clear answers to the minimum
achievable mean squared error (MSE), with the Cram\'er-Rao Bound
(CRB) being the most popular approach for benchmarking unbiased
estimators. The CRB is, typically, only valid conditional on having small
estimation errors. Although the CRB could be improved in the large-errors regime by the method of interval errors \cite{VanTrees07}, in general more sophisticated bounds should be explored in this region to handle the threshold phenomena (that is,
when the performance breaks down). Under such regime, Bayesian
bounds can be classified into pertaining to either the Ziv-Zakai or
the Weiss-Weinstein families \cite{VanTrees07}. However, all these
bounds assume that the model is perfectly specified, and thus they
provide a bound on the MSE performance under optimal conditions.

In the context of robust statistics, it has been said earlier that a
useful performance metric is to evaluate the efficiency of the
method by comparing its estimation performance, for instance
measured in terms of MSE, to the theoretical lower bound when the
model is perfectly specified. However, this bound might be too
optimistic, and thus unrealistic, under model mismatches. In
contrast, we are interested here in deriving estimation bounds
accounting for model departures.

In this paper we derive bounds on the achievable MSE in
the presence of model inaccuracies. The final goal being to compare
the performance of an estimation method (robust or non-robust) to
its theoretical bound when the model is wrongly specified. This
issue has been firstly addressed in \cite{Vuong86}, with an attempt
to generalize the CRB methodology. However, the results were not
conclusive and little research has been conducted in this important
direction. The approach in \cite{Vuong86} was interesting, but had
some limitations inherent to the CRB since inaccurate models
typically imply biased estimates and/or large-errors, and thus the
CRB may not be a valid approach in general. More recently, the idea
has been retaken in \cite{Richmond13,Fritsche2013}. Recently, in \cite{Richmond15}, the so called misspecified CRB (MCRB) was presented. The result has the feature of relating the bound with the Kullback-Leibler divergence between assumed and true distributions, giving an information-theoretic meaning to the bound. However, regardless its simple and convenient closed-form expression, a major drawback of the MCRB is the assumption of being able to compute the estimator's bias. Model misspecification typically causes biases in the estimation, and thus the standard CRB methodology has to be extended to account for such bias. This bias calculation could be particularly difficult in cases where, for instance, no closed-form expression exists for the estimator. Moreover, the bound is limited to the class of estimators sharing a specified bias. In this paper, we aim at obtaining bounds that are applicable to any type of estimators, whose only shared characteristic being that they are derived using an assumed model. We are interested in using the Ziv-Zakai bound (ZZB) methodology for this purpose. The ZZB has been previously considered to bound the MSE of estimators in misspecified models in \cite{Xu04}, thus implicitly taken bias into consideration. They restricted to the class of mismatch where the true distribution is wrongly parametrized. Parameters were divided into those being estimated and those assumed known. The latter being the possible cause for model inaccuracy. The authors derived a bound in this situation and applied it to the problem of direction finding using antenna arrays. In \cite{Gusi2014b}, a ZZB bound for model mismatching was derived for the particular case of time of arrival (TOA) estimation in the presence of unknown interference. The bound presented in the present article provides a more general result.

This paper derives a ZZB bound for the model mismatching problem
under additive Gaussian models. The derived bound is able to predict
the attainable performance of estimators under a number of model
inaccuracies. Particularly, these inaccuracies are considered in the
modeling of observations, that is the likelihood distribution of
data given the unknowns, while the a priori distribution is assumed
correctly elicited \cite{Rios00}. Namely, we have identified three
types of possible errors: $i$) when the specification of the
hyperparameters that are assumed known does not correspond to their
true values. Following the nomenclature in \cite{Xu04}, this case
corresponds to the background parameter mismatch; $ii$) when the,
possibly nonlinear, function relating the unknown parameters (those
we are bounding) and the observation differs from the actual
relationship; and $iii$) when the underlying noise distribution is
wrongly modeled, that is when the Gaussian distribution does not
reflect the true underlying law. The results are given first in its
general form and then particularized to special cases, where it is
easier to show closed-form analytical expressions and provide
insights on the results.

The remainder of the paper is organized as follows. Section
\ref{sec:problem} presents the problem and introduces some
mathematical notation. The more general result is given in Section
\ref{sec:ZZLB}, both for scalar and vector parameters. Then, the
special cases where the true noise distributions are Gaussian and
Gaussian mixture are discussed in Section \ref{sec:special_Gaussian}
and \ref{sec:special_Gaussian_mixture}, respectively. Section
\ref{sec:simus} illustrates the results with some examples and
Section \ref{sec:conc} concludes the paper with final remarks.

\section{Problem statement}
\label{sec:problem}

The additive Gaussian model is widely used in practical cases. In
this work we restrict to this family for the set of assumed models,
while the true underlying law being anything else. The unknown
parameter of interest, whose MSE we would like to bound, is denoted
by $\btheta \in  \Theta \subset \mathbb{R}^{n_\theta}$. The available
data is the random $K \times 1$ vector $\mathbf{x}$ which depends on
$\btheta$. Then, the assumed additive Gaussian model for
$\mathbf{x}$ is
\begin{equation}
\mathbf{x}=\mathbf{h}(\btheta)+\mathbf{n} \label{eq:esignal}
\end{equation}
where $\mathbf{h}(\cdot): \Theta \mapsto \mathbb{R}^K $ is a,
possibly nonlinear, function relating $\btheta$ to the observations.
$\mathbf{n}$ is a multivariate additive Gaussian noise term with
mean $\bm{\mu}$ and covariance matrix $\bm{\Sigma}$, conforming the
Gaussian hyperparameters $\bm{\Psi}=\{\bm{\mu} , \bm{\Sigma}\}$ that
can be either known or included in $\btheta$. The assumed
statistical model is therefore defined by the set of parameterized
distributions $\mathcal{M} = \{\mathbf{x} \sim p (\mathbf{x} |
\btheta) ~:~ \btheta \in \Theta \subset \mathbb{R}^{n_\theta}\}$,
where $p (\mathbf{x} | \btheta) =
\mathcal{N}(\mathbf{h}(\btheta);~\bm{\mu} , \bm{\Sigma})$ throughout
the paper.

The true data-generation process might differ from the accepted
description in (\ref{eq:esignal}). In the sequel, we use the
subscript $\ast$ to denote the hyperparameters, functions, or
distributions defining the true nature of the observations. Then,
the correctly specified model for $\mathbf{x}$ is defined by
\begin{equation}
\mathbf{x}=\mathbf{h}_\ast(\btheta)+\mathbf{n}_\ast \label{eq:signal}
\end{equation}
where $\mathbf{h}_\ast(\cdot): \Theta \mapsto \mathbb{R}^K $ is a,
possibly nonlinear, function of the unknown parameter $\btheta$, and
$\mathbf{n}_\ast$ is a random noise term with probability density
function (pdf) given by $p_{\ast}({\mathbf{n}_\ast})$ and with
hyperparameters gathered in $\bm{\Psi}_\ast$. Let us define the
statistical model under correct modeling assumptions as the set
$\mathcal{M}_\ast = \{\mathbf{x} \sim p_\ast (\mathbf{x} | \btheta)
~:~ \btheta \in \Theta \subset \mathbb{R}^{n_\theta}\}$. Notice
that, even when the distribution has been correctly selected, the
assumed model might have wrong hyperparameters $\bm{\Psi} \neq
\bm{\Psi}_\ast$. In case all (or some) of the parameters in
$\bm{\Psi}$ are not known, then they should be included in $\btheta$ and thus we end up in the classical estimation problem without mismatch (at least for the parameters in $\bm{\Psi}$ that are estimated).

When building estimators, an important component is the
log-likelihood function defined here as $\mathcal{L} (\mathbf{x} ;
\btheta, \mathcal{M}) = \log p (\mathbf{x} | \btheta)$ with
$\mathbf{x} \sim p_\ast (\mathbf{x} | \btheta)$ in reality. Then,
the maximum likelihood estimator (MLE) of $\btheta$ is
\begin{equation} \label{eq:quasiMLE}
  \hat{\btheta} = \arg \min_{\bm{\theta} \in \Theta} \mathcal{L} (\mathbf{x} ; \btheta, \mathcal{M})
\end{equation}
\noindent which is also referred to as the quasi-MLE when
$\mathcal{M} \neq \mathcal{M}_\ast$. It is known \cite{Huber67,Akaike73,
Sawa78, White82} that when the elements in $\mathbf{x}$ are
independent and identically distributed, then the quasi-MLE is also
the minimizer of the Kullback-Leibler Information Criterion (KLIC).
In other words, the quasi-MLE provides the smallest distance between
$\mathcal{M}$ and $\mathcal{M}_\ast$ in the sense that
\begin{equation}
  \hat{\btheta} \overset{a.s.}{\underset{K\rightarrow \infty}{\longrightarrow}} \arg \min_{\bm{\theta} \in \Theta} \mathbb{E}_\ast \left\{ \log \frac{p_\ast (\mathbf{x} | \btheta)}{p (\mathbf{x} | \btheta)}
  \right\} ~,
\end{equation}
\noindent where $\mathbb{E}_\ast \{\cdot\}$ denotes the expectation
with respect to the true distribution of the data, $p_\ast
(\mathbf{x} | \btheta)$. In this paper we aim at deriving estimation
bounds for estimators of the class in (\ref{eq:quasiMLE}).

\section{Ziv Zakai Lower Bound under model misspecification}
\label{sec:ZZLB} In this section we present a short review of the
ZZB. Unlike the CRB, the ZZB provides a bound on the MSE over the a
priori pdf of the unknown parameter. Moreover, the bound can
accomodate biased estimates and its application is therefore not
restricted to the family of unbiased estimators. The bound was first
derived in \cite{Chazan1975} for scalar parameters and subsequently
adapted to vector parameters in \cite{Bell1997}. The scalar and
multivariate versions of the classical bound are reviewed hereafter
in Sections \ref{sec:ZZLB:bound_scalar} and
\ref{sec:ZZLB:bound_vector}, respectively. Finally, we devote
Section \ref{sec:ZZLB:Pe} to evaluate the component of the bound
that incorporates the information regarding the model mismatching.
The latter providing the most general form of the bound for the
class of misspecified models described in Section \ref{sec:problem}.

\subsection{Fundamental estimation bound for scalar parameters}
\label{sec:ZZLB:bound_scalar}

Let us consider a scalar unknown parameter denoted by $\theta \in  \Theta \subset \mathbb{R}$. A lower-bound for the MSE of $\theta$
\begin{equation}
\mathbb{E}_\ast\{\epsilon^2\}=\mathbb{E}_\ast\left\{\left(\hat{\theta}-\theta\right)^2\right\}
\end{equation}
is envisaged. The ZZB can be obtained from the identity
\begin{equation}
\mathbb{E}_\ast\{\epsilon^2\}=\frac{1}{2}\int_{0}^{\infty}\mathbb{P}\left(\left|\epsilon\right|\geq\frac{h}{2}\right)h
~ dh
\end{equation}
and lower bounding
$\mathbb{P}\left(\left|\epsilon\right|\geq\frac{h}{2}\right)$, which is defined as
\begin{equation}
    \mathbb{P}\left(\left|\epsilon\right|\geq\frac{h}{2}\right) = \int \int \mathcal{I}\left(\left|\epsilon\right|\geq\frac{h}{2}\right) p_\ast(\mathbf{x},\theta) d\mathbf{x} d \theta
\end{equation}
\noindent in our problem. $\mathcal{I}(\cdot)$ denotes the indicator function of the argument even. The expression
$\mathbb{P}\left(\left|\epsilon\right|\geq\frac{h}{2}\right)$, with
$h\geq 0$, is related to a binary detection scheme with equally
probable hypotheses
\begin{IEEEeqnarray}{LL}
\mathcal{H}_1:\theta=\theta_o;&\mathbf{x}=\mathbf{h}(\theta)+\mathbf{n}|\theta=\theta_o\IEEEnonumber\\
\mathcal{H}_2:\theta=\theta_o+h;\:\:\:&\mathbf{x}=\mathbf{h}(\theta)+\mathbf{n}|\theta=\theta_o+h \label{eq:hyptest1}
\end{IEEEeqnarray}
when considering a suboptimal decision scheme, where the parameter
is first estimated and a nearest-neighbor decision is made
afterwards
\begin{equation}
\setlength{\nulldelimiterspace}{0pt} \hat{\mathcal{H}}=\left\{
\begin{IEEEeqnarraybox}[\relax][c]{l's}
\mathcal{H}_1,&if $\hat{\theta}\leq \theta_o+\frac{h}{2}$\\
\mathcal{H}_2,&if $\hat{\theta}> \theta_o+\frac{h}{2}$%
\end{IEEEeqnarraybox}\right.
\end{equation}

Intuitively, this test can be read as testing whether the estimator is performing good or committing some error larger than $h>0$, that is $\mathcal{H}_1$ or $\mathcal{H}_2$ respectively. The type of estimators we aim at bounding, $\hat{\theta}(\mathbf{x})$, use observations taken from the true model ($\mathbf{x} \sim p_\ast(\mathbf{x} | \theta)$). However, the estimator is derived considering the assumed model, $p(\mathbf{x} | \theta)$.
From a Bayesian theory perspective, the hypothesis testing in (\ref{eq:hyptest1}) is equivalent to the following model fitting problem
\begin{eqnarray}
\mathcal{H}_1: \mathbf{x} & \sim & p(\mathbf{x} | \theta_o) \label{eq:hyptest2}\\
\nonumber \mathcal{H}_2: \mathbf{x} & \sim & p(\mathbf{x} | \theta_o + h)
\end{eqnarray}
\noindent where the problem is to decide which distribution (i.e., $p(\mathbf{x} | \theta) $ or $p(\mathbf{x} | \theta + h) $) fits better the data $\mathbf{x}$ drawn from $p_\ast(\mathbf{x} | \theta)$. In general, the true distribution of measurements is not available. With these assumptions, the hypothesis test in (\ref{eq:hyptest2}) can be studied within the Bayesian statistics framework, where the possibility of performing hypothesis testing over data with unknown distribution is doable \cite[Chapter 6]{Bernardo09}. We can identify the problem as one about finding which parameters of a certain (assumed) distribution make the observations more probable, which is indeed what the estimator $\hat{\theta}(\mathbf{x})$ is aiming at.

Under this Bayesian perspective, and assuming uniform \emph{a priori} probabilities for the two hypotheses, the test in (\ref{eq:hyptest2}) -- consequently (\ref{eq:hyptest1}) -- can be solved by computing the likelihood ratio
\begin{equation}
\Lambda(\mathbf{x})=\frac{p(\mathbf{x}|\theta_o)}{p(\mathbf{x}|\theta_o+h)}\overset{\mathcal{H}_1}{\underset{\mathcal{H}_2}{\gtrless}}1
~, \label{eq:lrt_scalar}
\end{equation}
\noindent the evaluation of its minimum error probability $\mathbb{P}_e(a,a+h)$ being the objective of Section \ref{sec:ZZLB:Pe}.
The term $\mathbb{P}\left(\left|\epsilon\right|\geq\frac{h}{2}\right)$ can be
shown \cite{Bell1997} to be greater or equal to
\begin{equation}
\int_{-\infty}^{\infty}\left(p_{\theta}\left(\theta_o\right)+p_{\theta}\left(\theta_o+h\right)\right)\mathbb{P}_e(\theta_o,\theta_o+h)d\theta_o ~,
\end{equation}
where $p_{\theta}\left(\theta\right)$ is the a priori distribution
of the parameter of interest $\theta$. Assuming that
$p_{\theta}\left(\theta\right)$ follows a uniform distribution in
the interval $[0, T]$, the lower bound on the estimation error can
then be expressed as
\begin{equation}
\mathbb{E}\{\epsilon^2\}\geq\mbox{ZZB}=\frac{1}{T}\int_{0}^{T}h\int_{0}^{T-h}\mathbb{P}_e(\theta_o,\theta_o+h)d\theta_o
dh. \label{eq:ZZLB_scalar}
\end{equation}
Moreover, when $\mathbb{P}_e(\theta_o,\theta_o+h)$ is independent of
$\theta_o$ we can write $\mathbb{P}_e(h)$ instead. Under the latter
assumption, the ZZB reduces to evaluating the integral
\begin{equation}
\mbox{ZZB}=\frac{1}{T}\int_{0}^{T}h(T-h)\mathbb{P}_e(h)dh.
\label{eq:ZZLB_scalar_ind}
\end{equation}


\subsection{Fundamental estimation bound for vector parameters}
\label{sec:ZZLB:bound_vector} The bound provided above targets a
lower bound on the MSE for an unknown scalar parameter $\theta$
\cite{Bell1997}. Here we briefly present the derivation of the ZZB
for an unknown vector parameter $\btheta \in  \Theta \subset
\mathbb{R}^{n_\theta}$. For any estimator
$\hat{\boldsymbol{\theta}}$, the estimation error is given by
$\boldsymbol{\epsilon}=
\hat{\boldsymbol{\theta}}-\boldsymbol{\theta}$. We are interested in
a lower bound for the MSE of $\boldsymbol{\theta}$. The error
correlation matrix is obtained as
\begin{equation}
\mathbf{R}_{\boldsymbol{\epsilon}}=\mathbb{E}\left\{\left(\hat{\boldsymbol{\theta}}-\boldsymbol{\theta}\right)\left(\hat{\boldsymbol{\theta}}-\boldsymbol{\theta}\right)^\top\right\}=\mathbb{E}\left\{\boldsymbol{\epsilon}\boldsymbol{\epsilon}^\top\right\}~.
\end{equation}
A lower bound on
$\mathbf{a}^\top\mathbf{R}_{\boldsymbol{\epsilon}}\mathbf{a}$ is
envisaged for any $n_{\theta}$-dimensional vector  $\mathbf{a}$. If
a bound on a particular component of $\boldsymbol{\theta}$ is
required, $\mathbf{a}$ can be set to be a unit vector with a one in
the corresponding position and zeroes otherwise. The ZZB can be
obtained from the identity
\begin{equation}
\mathbf{a}^\top\mathbf{R}_{\boldsymbol{\epsilon}}\mathbf{a}=\frac{1}{2}\int_{0}^{\infty}\mathbb{P}\left(\left|\mathbf{a}^\top\epsilon\right|\geq\frac{h}{2}\right)hdh
\label{eq:zzlbound}
\end{equation}
and lower bounding
$\mathbb{P}\left(\left|\mathbf{a}^\top\boldsymbol{\epsilon}\right|\geq\frac{h}{2}\right)$.
The expression
$\mathbb{P}\left(\left|\mathbf{a}^\top\boldsymbol{\epsilon}\right|\geq\frac{h}{2}\right)$
is related to a binary detection scheme with equally probable
hypotheses
\begin{IEEEeqnarray}{LL}
\mathcal{H}_1:\boldsymbol{\theta}=\boldsymbol{\theta}_o;&\mathbf{x}=\mathbf{h}(\btheta)+\mathbf{n}|\boldsymbol{\theta}=\boldsymbol{\theta}_o\IEEEnonumber\\
\mathcal{H}_2:\boldsymbol{\theta}=\boldsymbol{\theta}_o+\boldsymbol{\delta};\:\:\:&\mathbf{x}=\mathbf{h}(\btheta)+\mathbf{n}|\boldsymbol{\theta}=\boldsymbol{\theta}_o+\boldsymbol{\delta}
\end{IEEEeqnarray}
with $\boldsymbol{\delta}$ satisfying
\begin{equation}
\mathbf{\mathbf{a}}^\top\boldsymbol{\delta}=h  \label{eq:delta}
\end{equation}
\noindent and $h$ non-negative. In fact,
$\mathbb{P}\left(\left|\mathbf{a}^\top\boldsymbol{\epsilon}\right|\geq\frac{h}{2}\right)$
can be lower bounded upon considering a suboptimal decision scheme
as in the scalar case, where the parameter is first estimated and a
nearest-neighbor decision is made afterwards
\begin{equation}
\setlength{\nulldelimiterspace}{0pt}
\hat{\mathcal{H}}=\left\{\begin{IEEEeqnarraybox}[\relax][c]{l's}
\mathcal{H}_1,&if $\mathbf{a}^\top\hat{\boldsymbol{\theta}}\leq \mathbf{a}^\top{\boldsymbol{\theta}_o}+\frac{h}{2}$\\
\mathcal{H}_2,&if $\mathbf{a}^\top\hat{\boldsymbol{\theta}}> \mathbf{a}^\top{\boldsymbol{\theta}_o}+\frac{h}{2}.$%
\end{IEEEeqnarraybox}\right.
\end{equation}
The probability of error for this suboptimum detector can be lower
bounded by the minimum error probability
$\mathbb{P}_e(\boldsymbol{\theta}_o,\boldsymbol{\theta}_o+\boldsymbol{\delta})$
given by the LRT
\begin{equation}
\Lambda(\mathbf{x})=\frac{p(\mathbf{x}|\boldsymbol{\theta}_o)}{p(\mathbf{x}|\boldsymbol{\theta}_o+\boldsymbol{\delta})}\overset{\mathcal{H}_1}{\underset{\mathcal{H}_2}{\gtrless}}1.
\label{eq:lrt}
\end{equation}

The term
$\mathbb{P}\left(\left|\mathbf{a}^\top\boldsymbol{\epsilon}\right|\geq\frac{h}{2}\right)$
can be shown \cite{Bell1997} to be greater or equal to
\begin{equation}
2\int\mbox{min}\left[p_{\boldsymbol{\theta}}\left(\boldsymbol{\theta}_o\right),p_{\boldsymbol{\theta}}\left(\boldsymbol{\theta}_o+\boldsymbol{\delta}\right)\right]\mathbb{P}_e(\boldsymbol{\theta}_o,\boldsymbol{\theta}_o+\boldsymbol{\delta})d\boldsymbol{\theta}_o,
\end{equation}
where $p_{\boldsymbol{\theta}}\left(\boldsymbol{\theta}\right)$ is,
again, the a priori distribution of $\boldsymbol{\theta}$. The bound
is valid for any $\boldsymbol{\delta}$ satisfying (\ref{eq:delta})
and the tightest one is found upon maximizing over all
$\boldsymbol{\delta}$ satisfying this constraint. This can be
written as
\begin{IEEEeqnarray}{LL}
\mathbb{P}\left(\left|\mathbf{a}^\top\boldsymbol{\epsilon}\right|\geq\frac{h}{2}\right)\geq\IEEEnonumber\\
\max_{\boldsymbol{\delta}:\mathbf{\mathbf{a}}^\top\boldsymbol{\delta}=h} 2\int\mbox{min}\left[p_{\boldsymbol{\theta}}\left(\boldsymbol{\theta}_o\right),p_{\boldsymbol{\theta}}\left(\boldsymbol{\theta}_o+\boldsymbol{\delta}\right)\right]\mathbb{P}_e(\boldsymbol{\theta}_o,\boldsymbol{\theta}_o+\boldsymbol{\delta})d\boldsymbol{\theta}_o.\IEEEnonumber\\
\label{eq:ZZLB_vector}
\end{IEEEeqnarray}
In some particular problems the minimum probability of error is not
a function of the parameter $\boldsymbol{\theta}_o$, but only of the
offset  between hypotheses. In this case, we have that
$\mathbb{P}_e(\boldsymbol{\theta}_o,\boldsymbol{\theta}_o+\boldsymbol{\delta})=\mathbb{P}_e(\boldsymbol{\delta})$
and the bound simplifies as
\begin{IEEEeqnarray}{LL}
\mathbb{P}\left(\left|\mathbf{a}^\top\boldsymbol{\epsilon}\right|\geq\frac{h}{2}\right)\geq
2
\max_{\boldsymbol{\delta}:\mathbf{\mathbf{a}}^\top\boldsymbol{\delta}=h}
A(\boldsymbol{\delta})\mathbb{P}_e(\boldsymbol{\delta})
\end{IEEEeqnarray}
where
\begin{equation}
A(\boldsymbol{\delta})=\int\mbox{min}\left[p_{\boldsymbol{\theta}}\left(\boldsymbol{\theta}_o\right),p_{\boldsymbol{\theta}}\left(\boldsymbol{\theta}_o+\boldsymbol{\delta}\right)\right]d\boldsymbol{\theta}_o.
\end{equation}
The bound for a probability of error independent of $\boldsymbol{\theta}_o$ is given by
\begin{equation}
\mathbf{a}^\top\mathbf{R}_{\boldsymbol{\epsilon}}\mathbf{a} \geq
\mbox{ZZB}=
\int_{0}^{\infty}\max_{\boldsymbol{\delta}:\mathbf{\mathbf{a}}^\top\boldsymbol{\delta}=h}
A(\boldsymbol{\delta})\mathbb{P}_e(\boldsymbol{\delta})h ~ dh.
\label{eq:ZZLB_vector_ind}
\end{equation}

\subsection{Minimum probability of error for misspecified models}
\label{sec:ZZLB:Pe} In this section we derive the minimum
probability of error associated with the rather general model
missmatch problem described in Section \ref{sec:problem}. For the
sake of generality, the derivation of the probability of error is
shown for the vector version of the bound, that is
$\mathbb{P}_e(\boldsymbol{\theta}_o,\boldsymbol{\theta}_o+\boldsymbol{\delta})$.
If one is interested in the scalar version of the bound, its
counterpart minimum error probability,
$\mathbb{P}_e(\theta_o,\theta_o+h)$, can be straightforwardly
obtained upon substituting the vector parameters
$\boldsymbol{\theta}_o$ and $\boldsymbol{\delta}$, with $\theta_o$
and $h$, respectively.

The minimum error probability
$\mathbb{P}_e(\boldsymbol{\theta}_o,\boldsymbol{\theta}_o+\boldsymbol{\delta})$
is given by the LRT in (\ref{eq:lrt}), which can be seen as a Bayesian classifier \cite[Chapter 2]{Theodoridis03}. For the sake of convenience,
the log-likelihood function definition in Section \ref{sec:problem}
is shortened as $\mathcal{L} ( \btheta) = \mathcal{L} (\mathbf{x} ;
\btheta, \mathcal{M})$, where the dependence with observations
$\mathbf{x} \sim p_\ast (\mathbf{x} | \btheta)$ and the assumed
model $\mathcal{M}$ is omitted. The log-likelihood ratio (LLR) can
be obtained upon taking the logarithm
\begin{IEEEeqnarray}{LLL}
\ln\Lambda(\mathbf{x})&=&\ln{p(\mathbf{x}|\boldsymbol{\theta}_o)}-\ln{p(\mathbf{x}|\boldsymbol{\theta}_o+\boldsymbol{\delta})}\IEEEnonumber\\
&=&\mathcal{L}(\boldsymbol{\theta}_o)-\mathcal{L}(\boldsymbol{\theta}_o+\boldsymbol{\delta})
\overset{\mathcal{H}_1}{\underset{\mathcal{H}_2}{\gtrless}}0.
\label{eq:lnlrt}
\end{IEEEeqnarray}

The log-likelihood function of a multivariate normal distribution,
neglecting the irrelevant constant terms, is given by
\begin{equation}
\mathcal{L}(\btheta)=-\frac{1}{2}\left(\mathbf{x}-\mathbf{h}(\btheta)-\boldsymbol{\mu}\right)^\top\boldsymbol{\Sigma}^{-1}\left(\mathbf{x}-\mathbf{h}(\btheta)-\boldsymbol{\mu}\right)
\label{eq:Likelihood}
\end{equation}
where $\mathbf{x}$ is drawn from the true law $\mathbf{x} \sim
p_\ast (\mathbf{x} | \btheta)$ defined in (\ref{eq:signal}). The
minimum error probability is then
\begin{IEEEeqnarray}{LLL}
\mathbb{P}_e(\boldsymbol{\theta}_o,\boldsymbol{\theta}_o+\boldsymbol{\delta})&=&\mathbb{P}\left(\ln\Lambda(\mathbf{x})<0|\mathcal{H}_1\right)\mathbb{P}\left(\mathcal{H}_1\right)\IEEEnonumber\\
&&+\mathbb{P}\left(\ln\Lambda(\mathbf{x})>0|\mathcal{H}_2\right)\mathbb{P}\left(\mathcal{H}_2\right)\IEEEnonumber\\
&=&\frac{1}{2}\mathbb{P}\left(\ln\Lambda(\mathbf{x})<0|\mathcal{H}_1\right)\IEEEnonumber\\
&&+\frac{1}{2}\mathbb{P}\left(\ln\Lambda(\mathbf{x})>0|\mathcal{H}_2\right),\label{eq:Pe_general}
\end{IEEEeqnarray}
where equally likely hypotheses are assumed for the second equality.
The remaining probabilities can be obtained as
\begin{IEEEeqnarray}{ll}
\mathbb{P}\left(\ln\Lambda(\mathbf{x})<0|\mathcal{H}_1\right)&=\mathbb{P}\left(\mathcal{L}(\boldsymbol{\theta}_o)-\mathcal{L}(\boldsymbol{\theta}_o+\boldsymbol{\delta})<0|\btheta=\boldsymbol{\theta}_o\right)\label{eq:P10}\\
\mathbb{P}\left(\ln\Lambda(\mathbf{x})>0|\mathcal{H}_2\right)&=\mathbb{P}\left(\mathcal{L}(\boldsymbol{\theta}_o)-\mathcal{L}(\boldsymbol{\theta}_o+\boldsymbol{\delta})>0|\btheta=\boldsymbol{\theta}_o+\boldsymbol{\delta}\right).\IEEEnonumber\label{eq:P01}\\
\end{IEEEeqnarray}

The log-likelihood function evaluated at $\boldsymbol{\theta}_o$ in
(\ref{eq:Likelihood}) can be expanded as
\begin{IEEEeqnarray}{LLL}
\mathcal{L}(\boldsymbol{\theta}_o)&=&-\frac{1}{2}\left(\mathbf{x}-\mathbf{h}(\boldsymbol{\theta}_o)-\boldsymbol{\mu}\right)^\top\boldsymbol{\Sigma}^{-1}\left(\mathbf{x}-\mathbf{h}(\boldsymbol{\theta}_o)-\boldsymbol{\mu}\right)\IEEEnonumber\\
&=&-\frac{1}{2}\left(\mathbf{x}^\top\boldsymbol{\Sigma}^{-1}\mathbf{x}+(\mathbf{h}(\boldsymbol{\theta}_o)+\boldsymbol{\mu})^\top\boldsymbol{\Sigma}^{-1}\left(\mathbf{h}(\boldsymbol{\theta}_o)+\boldsymbol{\mu}\right)\right)\IEEEnonumber\\
&&+\mathbf{x}^\top\boldsymbol{\Sigma}^{-1}\left(\mathbf{h}(\boldsymbol{\theta}_o)+\boldsymbol{\mu}\right)
~, \label{eq:Deltaa}
\end{IEEEeqnarray}
where
\begin{IEEEeqnarray}{LLL}
\mathbf{x}^\top\boldsymbol{\Sigma}^{-1}(\mathbf{h}(\boldsymbol{\theta}_o)+\boldsymbol{\mu})&=&\left(\mathbf{h}_\ast(\btheta)+\mathbf{n}_\ast\right)^\top\boldsymbol{\Sigma}^{-1}(\mathbf{h}(\boldsymbol{\theta}_o)+\boldsymbol{\mu})\IEEEnonumber\\
&=&\mathbf{h}_\ast^\top(\btheta)\boldsymbol{\Sigma}^{-1}(\mathbf{h}(\boldsymbol{\theta}_o)+\boldsymbol{\mu})\IEEEnonumber\\
&&+\mathbf{n}_\ast^\top\boldsymbol{\Sigma}^{-1}(\mathbf{h}(\boldsymbol{\theta}_o)+\boldsymbol{\mu})
~.
\end{IEEEeqnarray}
\noindent Notice that the observation vector is substituted by its
model under $\mathcal{M}_\ast$, that is
$\mathbf{x}=\mathbf{h}_\ast(\btheta)+\mathbf{n}_\ast$, where the
parameter $\theta$ is left arbitrary and will be particularized
later depending on the hypothesis of the LRT.

Similarly, the log-likelihood function evaluated at
$\boldsymbol{\theta}_o+\boldsymbol{\delta}$ is
\begin{IEEEeqnarray}{LLL}
\mathcal{L}(\boldsymbol{\theta}_o+\boldsymbol{\delta})&=&-\frac{1}{2}\mathbf{x}^\top\boldsymbol{\Sigma}^{-1}\mathbf{x}\IEEEnonumber\\
&&-\frac{1}{2}(\mathbf{h}(\boldsymbol{\theta}_o+\boldsymbol{\delta})+\boldsymbol{\mu})^\top\boldsymbol{\Sigma}^{-1}(\mathbf{h}(\boldsymbol{\theta}_o+\boldsymbol{\delta})+\boldsymbol{\mu})\IEEEnonumber\\
&&+\mathbf{x}^\top\boldsymbol{\Sigma}^{-1}(\mathbf{h}(\boldsymbol{\theta}_o+\boldsymbol{\delta})+\boldsymbol{\mu})~,
\label{eq:Deltaah}
\end{IEEEeqnarray}
where
\begin{IEEEeqnarray}{LLL}
\mathbf{x}^\top\boldsymbol{\Sigma}^{-1}\mathbf{h}(\boldsymbol{\theta}_o+\boldsymbol{\delta})&=&\mathbf{h}_\ast^\top(\btheta)\boldsymbol{\Sigma}^{-1}(\mathbf{h}(\boldsymbol{\theta}_o+\boldsymbol{\delta})+\boldsymbol{\mu})\IEEEnonumber\\
&&+\mathbf{n}_\ast^\top\boldsymbol{\Sigma}^{-1}(\mathbf{h}(\boldsymbol{\theta}_o+\boldsymbol{\delta})+\boldsymbol{\mu})
~.
\end{IEEEeqnarray}

From (\ref{eq:Deltaa}) and (\ref{eq:Deltaah}) we can compute $\ln
\Lambda(\mathbf{x})=\mathcal{L}(\boldsymbol{\theta}_o)-\mathcal{L}(\boldsymbol{\theta}_o+\boldsymbol{\delta})$
as
\begin{IEEEeqnarray}{LLL}
\ln\Lambda(\mathbf{x})&=&\frac{1}{2}(\mathbf{h}(\boldsymbol{\theta}_o+\boldsymbol{\delta})+\boldsymbol{\mu})^\top\boldsymbol{\Sigma}^{-1}(\mathbf{h}(\boldsymbol{\theta}_o+\boldsymbol{\delta})+\boldsymbol{\mu})\IEEEnonumber\\
&&-\frac{1}{2}(\mathbf{h}(\boldsymbol{\theta}_o)+\boldsymbol{\mu})^\top\boldsymbol{\Sigma}^{-1}(\mathbf{h}(\boldsymbol{\theta}_o)+\boldsymbol{\mu})\IEEEnonumber\\
&&+\mathbf{h}_\ast^\top(\btheta)\boldsymbol{\Sigma}^{-1}(\mathbf{h}(\boldsymbol{\theta}_o)-\mathbf{h}(\boldsymbol{\theta}_o+\boldsymbol{\delta}))\IEEEnonumber\\
&&+\mathbf{n}_\ast^\top\boldsymbol{\Sigma}^{-1}(\mathbf{h}(\boldsymbol{\theta}_o)-\mathbf{h}(\boldsymbol{\theta}_o+\boldsymbol{\delta}))
\end{IEEEeqnarray}

The probability in (\ref{eq:P10}) yields
\begin{IEEEeqnarray}{LLL}
\mathbb{P}\left(\ln\Lambda(\mathbf{x})<0|\mathcal{H}_1\right)&=&\mathbb{P}\left(S(\boldsymbol{\theta}_o,\boldsymbol{\delta})+n<0\right)
\label{eq:Ps1n}
\end{IEEEeqnarray}
where
\begin{IEEEeqnarray}{LLL}
S(\btheta,\boldsymbol{\delta})&=&\frac{1}{2}(\mathbf{h}(\boldsymbol{\theta}_o+\boldsymbol{\delta})+\boldsymbol{\mu})^\top\boldsymbol{\Sigma}^{-1}(\mathbf{h}(\boldsymbol{\theta}_o+\boldsymbol{\delta})+\boldsymbol{\mu})\IEEEnonumber\\
&&-\frac{1}{2}(\mathbf{h}(\boldsymbol{\theta}_o)+\boldsymbol{\mu})^\top\boldsymbol{\Sigma}^{-1}(\mathbf{h}(\boldsymbol{\theta}_o)+\boldsymbol{\mu})\IEEEnonumber\\
&&+\mathbf{h}_\ast^\top(\btheta)\boldsymbol{\Sigma}^{-1}(\mathbf{h}(\boldsymbol{\theta}_o)-\mathbf{h}(\boldsymbol{\theta}_o+\boldsymbol{\delta}))
\label{eq:S}
\end{IEEEeqnarray}
and
\begin{equation}\label{eq:n_add}
n=\mathbf{n}_\ast^\top\boldsymbol{\Sigma}^{-1}(\mathbf{h}(\boldsymbol{\theta}_o)-\mathbf{h}(\boldsymbol{\theta}_o+\boldsymbol{\delta})).
\end{equation}
with pdf $p_n(n)$. The pdf of $n$ is in general difficult to obtain,
as $n$ is given by the sum of noise samples from $\mathbf{n}_\ast
\sim p_\ast (\mathbf{n}_\ast)$ multiplied by the elements in
$\boldsymbol{\Sigma}^{-1}(\mathbf{h}(\boldsymbol{\theta}_o)-\mathbf{h}(\boldsymbol{\theta}_o+\boldsymbol{\delta}))$.
Moreover, $p_n(n)$ needs to be generated for different values of
$\boldsymbol{\theta}_o$ and $\boldsymbol{\delta}$, as the bound
integrates over these parameters. For arbitrary noise distributions,
the use of Monte Carlo methods would be recommended. The probability
in (\ref{eq:Ps1n}) can be computed as
\begin{IEEEeqnarray}{LLL}
\mathbb{P}\left(S(\boldsymbol{\theta}_o,\boldsymbol{\delta})+n<0\right)&=&\int_{-\infty}^{-S(\boldsymbol{\theta}_o,\boldsymbol{\delta})}p_n(x)dx
\label{eq:ints1}
\end{IEEEeqnarray}
In the same way, the probability in (\ref{eq:P01}) leads to
\begin{IEEEeqnarray}{LLL}
\mathbb{P}\left(\ln\Lambda(\mathbf{x})>0|\mathcal{H}_2\right)&=&\mathbb{P}\left(-S(\boldsymbol{\theta}_o+\boldsymbol{\delta},\boldsymbol{\delta})-n<0\right) \IEEEnonumber \\ 
{}&=&
\int_{-S(\boldsymbol{\theta}_o+\boldsymbol{\delta},\boldsymbol{\delta})}^{\infty}p_n(x)dx~,
\label{eq:ints2}
\end{IEEEeqnarray}
and, from substitution of (\ref{eq:ints1}) and (\ref{eq:ints2}) in
(\ref{eq:Pe_general}), we have that the minimum probability of error
is equal to
\begin{IEEEeqnarray}{LLL}
\mathbb{P}_e(\boldsymbol{\theta}_o,\boldsymbol{\theta}_o+\boldsymbol{\delta})&=&\frac{1}{2}\int_{-\infty}^{-S(\boldsymbol{\theta}_o,\boldsymbol{\delta})}p_n(x)dx\IEEEnonumber\\
&+&\frac{1}{2}\int_{-S(\boldsymbol{\theta}_o+\boldsymbol{\delta},\boldsymbol{\delta})}^{\infty}p_n(x)dx
\label{eq:P_e_vector}
\end{IEEEeqnarray}
Upon inserting the above result into (\ref{eq:ZZLB_vector}) the
lower bound for the general case is derived. Summarizing, this
general bound deals with misspecified models by means of different
noise distributions, and different functions of $\btheta$ (possibly
nonlinear). Regarding the noise, the assumed  model $\mathbf{n}$ is
set to a multivariate additive Gaussian distribution. This
consideration impacts in (\ref{eq:Likelihood}), where the
log-likelihood of a Gaussian multivariate distribution is used. The
correctly specified noise $\mathbf{n}_\ast$ is set to have arbitrary
pdf $p_{\ast}(\mathbf{n}_{\ast})$. The impact of the true noise
distribution $p_{\ast}(\mathbf{n}_{\ast})$ is modeled by means of
the scalar variable $n$ in (\ref{eq:n_add}), which is a
transformation of the original vector random variable
$\mathbf{n}_\ast$. Last, the correctly specified function of
$\btheta$, $\mathbf{h}_\ast(\cdot)$, appears in the last terms of
expression $S(\boldsymbol{\theta},\boldsymbol{\delta})$ in
(\ref{eq:S}).

\section{Special Cases: noise follows a Gaussian distribution} \label{sec:special_Gaussian}

In this section we particularize the minimum error probability result
from (\ref{eq:P_e_vector}) when $\mathbf{n}_\ast$ is a multivariate
additive Gaussian process. This consideration simplifies the
computation of the pdf $p_n(n)$ of $n$ in (\ref{eq:n_add}), since
the sum of independent random variables that are normally
distributed is also normally distributed.

Let the correctly specified noise signal $\mathbf{n}_\ast$ be a
multivariate additive Gaussian noise term with mean
$\boldsymbol{\mu}_\ast$ and covariance matrix
$\boldsymbol{\Sigma}_\ast$. In general, the assumed noise parameters
$\boldsymbol{\mu}$ and $\boldsymbol{\Sigma}$ might differ from
$\boldsymbol{\mu}_\ast$ and $\boldsymbol{\Sigma}_\ast$. The scalar
random variable $n$ in (\ref{eq:n_add}) is then normally
distributed, with mean $\mu_n$
\begin{IEEEeqnarray}{LLL}
\mu_n=\mathbb{E}\{n\}=\boldsymbol{\mu}_\ast^\top\boldsymbol{\Sigma}^{-1}(\mathbf{h}(\boldsymbol{\btheta}_o)-\mathbf{h}(\boldsymbol{\btheta}_o+\boldsymbol{\delta})).
\label{eq:meann}
\end{IEEEeqnarray}
and variance $\sigma_n^2$ equal to
\begin{IEEEeqnarray}{L}
\sigma_n^2=\mathbb{E}\{(n-\mu_n)^2\}\IEEEnonumber\\
=(\mathbf{h}(\boldsymbol{\btheta}_o)-\mathbf{h}(\boldsymbol{\btheta}_o+\boldsymbol{\delta}))^\top\boldsymbol{\Sigma}^{-1}\boldsymbol{\Sigma}_\ast\boldsymbol{\Sigma}^{-1}(\mathbf{h}(\boldsymbol{\btheta}_o)-\mathbf{h}(\boldsymbol{\btheta}_o+\boldsymbol{\delta})),\IEEEnonumber\\
\label{eq:variancen}
\end{IEEEeqnarray}
where
\begin{IEEEeqnarray}{LLL}
\boldsymbol{\Sigma}_\ast&=&\mathbb{E}\{(\mathbf{n}_\ast-\boldsymbol{\mu}_\ast)(\mathbf{n}_\ast-\boldsymbol{\mu}_\ast)^\top\}.
\end{IEEEeqnarray}
Let us define
\begin{IEEEeqnarray}{LLL}
Z(\boldsymbol{\btheta},\boldsymbol{\delta})&=&S(\boldsymbol{\btheta},\boldsymbol{\delta})+\mu_n\IEEEnonumber\\
&=&\frac{1}{2}(\mathbf{h}(\boldsymbol{\theta}_o+\boldsymbol{\delta})+\boldsymbol{\mu})^\top\boldsymbol{\Sigma}^{-1}(\mathbf{h}(\boldsymbol{\theta}_o+\boldsymbol{\delta})+\boldsymbol{\mu})\IEEEnonumber\\
&&-\frac{1}{2}(\mathbf{h}(\boldsymbol{\theta}_o)+\boldsymbol{\mu})^\top\boldsymbol{\Sigma}^{-1}(\mathbf{h}(\boldsymbol{\theta}_o)+\boldsymbol{\mu})\IEEEnonumber\\
&&+(\mathbf{h}_\ast(\boldsymbol{\theta})+\boldsymbol{\mu}_\ast)^\top\boldsymbol{\Sigma}^{-1}(\mathbf{h}(\boldsymbol{\theta}_o)-\mathbf{h}(\boldsymbol{\theta}_o+\boldsymbol{\delta})).
\label{eq:Zbthetadelta}
\end{IEEEeqnarray}

With this definition, the probabilities in (\ref{eq:Ps1n}) and
(\ref{eq:ints2}) can be computed as
\begin{IEEEeqnarray}{LLL}
\mbox{Pr}\left(S(\boldsymbol{\btheta}_o,\boldsymbol{\delta})+n<0\right)&=&Q\left(\frac{Z(\boldsymbol{\btheta}_o,\boldsymbol{\delta})}{\sigma_n}\right)
\label{eq:Qs1}
\end{IEEEeqnarray}
and
\begin{IEEEeqnarray}{LLL}
\mbox{Pr}\left(-S(\boldsymbol{\btheta}_o+\boldsymbol{\delta},\boldsymbol{\delta})-n<0\right)&=&Q\left(\frac{-Z(\boldsymbol{\btheta}_o+\boldsymbol{\delta},\boldsymbol{\delta})}{\sigma_n}\right),
\label{eq:Qs2}
\end{IEEEeqnarray}
where $Q(x)=(1/\sqrt{2\pi})\int_{x}^{\infty}\exp(-t^2/2))dt$ is the
Q-function, expressed in terms of the complementary error function
as $Q(x)=(1/2)\mbox{erfc}(x/\sqrt{2})$. The minimum probability of
error is then
\begin{IEEEeqnarray}{LLL}
P_{e}(\boldsymbol{\btheta}_o,\boldsymbol{\btheta}_o+\boldsymbol{\delta})&=&\frac{1}{2}Q\left(\frac{Z(\boldsymbol{\btheta}_o,\boldsymbol{\delta})}{\sigma_n}\right)\IEEEnonumber\\
&&+\frac{1}{2}Q\left(\frac{-Z(\boldsymbol{\btheta}_o+\boldsymbol{\delta},\boldsymbol{\delta})}{\sigma_n}\right)
\label{eq:Pe_gaussian}
\end{IEEEeqnarray}
and the bound is found by replacing the expression of $P_{e}(\boldsymbol{\btheta}_o,\boldsymbol{\btheta}_o+\boldsymbol{\delta})$
in (\ref{eq:ZZLB_vector}) with the above result.

In the following subsections more specific assumptions are taken
into account. First, we consider the case where the functions
$\mathbf{h}(\cdot)$ and  $\mathbf{h}_\ast(\cdot)$ are equal. This is the case of background parameter mismatch. Second, we treat the case where functions $\mathbf{h}(\cdot)$ and $\mathbf{h}_\ast(\cdot)$ are actually linear functions of $\btheta$. This case corresponds to both functional and background parameter
mismatches. And last, the case where the functions are both equal and linear, where we have again only background parameter mismatch. The latter case is then also particularized for a scenario where the noise mean of both the assumed and the correctly specified
models match, and for a scenario where the noise mean and
covariance matrix of both models coincide. This is the classical
result under perfectly-matched models \cite{Kay:1993:FSS:151045}, which we
include here as a sanity check of the more general bound under model
discrepancies.

\subsection{Functions of $\btheta$ are equal}
Let us assume now that the assumed functions of the unknown
parameter $\btheta$ coincide with the true transformation. We have
then that
\begin{IEEEeqnarray}{LLL}
\mathbf{h}(\btheta)&=&\mathbf{h}_\ast(\btheta)
\end{IEEEeqnarray}
The mean and the covariance of $n$ are now equal to
\begin{IEEEeqnarray}{LLL}
\mu_n=\mathbb{E}\{n\}=\boldsymbol{\mu}_\ast^\top\boldsymbol{\Sigma}^{-1}(\mathbf{h}_\ast(\boldsymbol{\btheta}_o)-\mathbf{h}_\ast(\boldsymbol{\btheta}_o+\boldsymbol{\delta})),
\end{IEEEeqnarray}
and
\begin{IEEEeqnarray}{LLL}
\sigma_n^2&=&(\mathbf{h}_\ast(\boldsymbol{\btheta}_o)-\mathbf{h}_\ast(\boldsymbol{\btheta}_o+\boldsymbol{\delta}))^\top\boldsymbol{\Sigma}^{-1}\boldsymbol{\Sigma}_\ast\IEEEnonumber\\
&&\cdot\boldsymbol{\Sigma}^{-1}(\mathbf{h}_\ast(\boldsymbol{\btheta}_o)-\mathbf{h}_\ast(\boldsymbol{\btheta}_o+\boldsymbol{\delta}))~,
\label{eq:variancen_equal}
\end{IEEEeqnarray}
\noindent respectively. The functions
$Z(\boldsymbol{\btheta}_o,\boldsymbol{\delta})$ and
$Z(\boldsymbol{\btheta}_o+\boldsymbol{\delta},\boldsymbol{\delta})$
can be reexpressed after some algebraic manipulation as
\begin{IEEEeqnarray}{lll}
Z(\boldsymbol{\btheta}_o,\boldsymbol{\delta})&=&\frac{1}{2}(\mathbf{h}_\ast(\boldsymbol{\theta}_o+\boldsymbol{\delta})-\mathbf{h}_\ast(\boldsymbol{\theta}_o))^\top\boldsymbol{\Sigma}^{-1}\IEEEnonumber\\
&&\cdot(\mathbf{h}_\ast(\boldsymbol{\theta}_o+\boldsymbol{\delta})-\mathbf{h}_\ast(\boldsymbol{\theta}_o))\IEEEnonumber\\
&&+(\mathbf{h}_\ast(\boldsymbol{\theta}_o+\boldsymbol{\delta})-\mathbf{h}_\ast(\boldsymbol{\theta}_o))^\top\boldsymbol{\Sigma}^{-1}(\boldsymbol{\mu}-\boldsymbol{\mu}_\ast)
\label{eq:Z1_equal}
\end{IEEEeqnarray}
and
\begin{IEEEeqnarray}{LLL}
Z(\boldsymbol{\btheta}_o+\boldsymbol{\delta},\boldsymbol{\delta})&=&-\frac{1}{2}(\mathbf{h}_\ast(\boldsymbol{\theta}_o+\boldsymbol{\delta})-\mathbf{h}_\ast(\boldsymbol{\theta}_o))^\top\boldsymbol{\Sigma}^{-1}\IEEEnonumber\\
&&\cdot(\mathbf{h}_\ast(\boldsymbol{\theta}_o+\boldsymbol{\delta})-\mathbf{h}_\ast(\boldsymbol{\theta}_o))\IEEEnonumber\\
&&+(\mathbf{h}_\ast(\boldsymbol{\theta}_o+\boldsymbol{\delta})-\mathbf{h}_\ast(\boldsymbol{\theta}_o))^\top\boldsymbol{\Sigma}^{-1}(\boldsymbol{\mu}-\boldsymbol{\mu}_\ast)\IEEEnonumber\\
\label{eq:Z2_equal}
\end{IEEEeqnarray}
The probability of error is found as in (\ref{eq:Pe_gaussian}) but
using the evaluated expressions of
$Z(\boldsymbol{\btheta},\boldsymbol{\delta})$ shown in
(\ref{eq:Z1_equal}) and (\ref{eq:Z2_equal}) and the variance in
(\ref{eq:variancen_equal}). The bound is obtained after inserting
this error probability in (\ref{eq:ZZLB_vector}).

\subsection{Functions of $\btheta$ are linear}
Let us assume now that the functions of the unknown parameter
$\btheta$ are linear. We have then that
\begin{IEEEeqnarray}{LLL}
\mathbf{h}(\btheta)&=&\mathbf{H}\btheta\IEEEnonumber\\
\mathbf{h}_\ast(\btheta)&=&\mathbf{H}_\ast\btheta
\end{IEEEeqnarray}
where $\mathbf{H}$ and $\mathbf{H}_\ast \in \mathbb{R}^{n_\theta
\times n_\theta}$ are both matrices denoting the assumed and
correctly specified linear functions of $\btheta$, respectively. The
function $S(\boldsymbol{\btheta},\boldsymbol{\delta})$ from
(\ref{eq:S}) yields to
\begin{IEEEeqnarray}{LLL}
S(\boldsymbol{\btheta},\boldsymbol{\delta})&=&\frac{1}{2}\boldsymbol{\delta}^\top\mathbf{H}^\top\boldsymbol{\Sigma}^{-1}\mathbf{H}\boldsymbol{\delta}+\boldsymbol{\delta}^\top\mathbf{H}^\top\boldsymbol{\Sigma}^{-1}(\mathbf{H}\boldsymbol{\btheta}_o+\boldsymbol{\mu})\IEEEnonumber\\
&&-\boldsymbol{\btheta}^\top\mathbf{H}_\ast^\top\boldsymbol{\Sigma}^{-1}\mathbf{H}\boldsymbol{\delta}.
\end{IEEEeqnarray}
The mean from (\ref{eq:meann}) is equal to
\begin{IEEEeqnarray}{LLL}
\mu_n=-\boldsymbol{\mu}_\ast^\top\boldsymbol{\Sigma}^{-1}\mathbf{H}\boldsymbol{\delta}=-\boldsymbol{\delta}^\top\mathbf{H}^\top\boldsymbol{\Sigma}^{-1}\boldsymbol{\mu}_\ast,
\end{IEEEeqnarray}
and the variance in (\ref{eq:variancen}) is now equal to
\begin{IEEEeqnarray}{LLL}
\sigma_n^2&=&\boldsymbol{\delta}^\top\mathbf{H}^\top\boldsymbol{\Sigma}^{-1}\boldsymbol{\Sigma}_\ast\boldsymbol{\Sigma}^{-1}\mathbf{H}\boldsymbol{\delta}.
\label{eq:variancen_linear}
\end{IEEEeqnarray}
Note that now the expression $Z(\boldsymbol{\btheta},\boldsymbol{\delta})=S(\boldsymbol{\btheta},\boldsymbol{\delta})+\mu_n$ can be
written as a function of the difference between the means of the
assumed and correct noise $\boldsymbol{\mu}-\boldsymbol{\mu}_\ast$
\begin{IEEEeqnarray}{LLL}
Z(\boldsymbol{\btheta},\boldsymbol{\delta})&=&\frac{1}{2}\boldsymbol{\delta}^\top\mathbf{H}^\top\boldsymbol{\Sigma}^{-1}\mathbf{H}\boldsymbol{\delta}+\boldsymbol{\delta}^\top\mathbf{H}^\top\boldsymbol{\Sigma}^{-1}(\mathbf{H}\boldsymbol{\btheta}_o+\boldsymbol{\mu}-\boldsymbol{\mu}_\ast)\IEEEnonumber\\
&&-\boldsymbol{\btheta}^\top\mathbf{H}_\ast^\top\boldsymbol{\Sigma}^{-1}\mathbf{H}\boldsymbol{\delta}.
\label{eq:Z1_linear}
\end{IEEEeqnarray}
Upon inserting the expression from (\ref{eq:Z1_linear}) and the
variance from (\ref{eq:variancen_linear}) in (\ref{eq:Pe_gaussian})
the probability of error is found for the linear case.

\subsection{Functions of $\btheta$ are equal and linear}
\label{functions_equal_linear}
Let us now suppose that the linear function of $\btheta$ from the assumed model
$\mathbf{H}$ matches with the correctly specified model $\mathbf{H}_\ast$. So, we
have that
\begin{equation}
\mathbf{H}=\mathbf{H}_\ast
\end{equation}
Function $Z(\boldsymbol{\btheta}_o,\boldsymbol{\delta})$  can be written as
\begin{IEEEeqnarray}{LLL}
Z(\boldsymbol{\btheta}_o,\boldsymbol{\delta})&=&\frac{1}{2}\boldsymbol{\delta}^\top\mathbf{H}_{\ast}^\top\boldsymbol{\Sigma}^{-1}\mathbf{H}_{\ast}\boldsymbol{\delta}+\boldsymbol{\delta}^\top\mathbf{H}_{\ast}^\top\boldsymbol{\Sigma}^{-1}(\mathbf{H}_{\ast}\boldsymbol{\btheta}_o+\boldsymbol{\mu}-\boldsymbol{\mu}_\ast)\IEEEnonumber\\
&&-\boldsymbol{\btheta}_o^\top\mathbf{H}_\ast^\top\boldsymbol{\Sigma}^{-1}\mathbf{H}_{\ast}\boldsymbol{\delta}\IEEEnonumber\\
&=&\frac{1}{2}\boldsymbol{\delta}^\top\mathbf{H}_{\ast}^\top\boldsymbol{\Sigma}^{-1}\mathbf{H}_{\ast}\boldsymbol{\delta}+\boldsymbol{\delta}^\top\mathbf{H}_{\ast}^\top\boldsymbol{\Sigma}^{-1}(\boldsymbol{\mu}-\boldsymbol{\mu}_\ast)\IEEEnonumber\\
&=&Z(\boldsymbol{\delta}).
\end{IEEEeqnarray}
Note that the dependency on $\boldsymbol{\btheta}_o$ has been
removed since the outcome does no depend on $\boldsymbol{\btheta}_o$
anymore. Similarly, function
$Z(\boldsymbol{\btheta}_o+\boldsymbol{\delta},\boldsymbol{\delta})$
is given by
\begin{IEEEeqnarray}{LLL}
Z(\boldsymbol{\btheta}_o+\boldsymbol{\delta},\boldsymbol{\delta})&=&\frac{1}{2}\boldsymbol{\delta}^\top\mathbf{H}_{\ast}^\top\boldsymbol{\Sigma}^{-1}\mathbf{H}_{\ast}\boldsymbol{\delta}+\boldsymbol{\delta}^\top\mathbf{H}_{\ast}^\top\boldsymbol{\Sigma}^{-1}\IEEEnonumber\\
&&\cdot(\mathbf{H}_{\ast}\boldsymbol{\btheta}_o+\boldsymbol{\mu}-\boldsymbol{\mu}_\ast)\IEEEnonumber\\
&&-(\boldsymbol{\btheta}_o+\boldsymbol{\delta})^\top\mathbf{H}_\ast^\top\boldsymbol{\Sigma}^{-1}\mathbf{H}_{\ast}\boldsymbol{\delta}\IEEEnonumber\\
&=&-\frac{1}{2}\boldsymbol{\delta}^\top\mathbf{H}_{\ast}^\top\boldsymbol{\Sigma}^{-1}\mathbf{H}_{\ast}\boldsymbol{\delta}+\boldsymbol{\delta}^\top\mathbf{H}_{\ast}^\top\boldsymbol{\Sigma}^{-1}(\boldsymbol{\mu}-\boldsymbol{\mu}_\ast)\IEEEnonumber\\
&=&-Z(\boldsymbol{-\delta}).
\end{IEEEeqnarray}
The mean and the variance are now equal to
\begin{IEEEeqnarray}{LLL}
\mu_n=-\boldsymbol{\mu}_\ast^\top\boldsymbol{\Sigma}^{-1}\mathbf{H}_\ast\boldsymbol{\delta},
\end{IEEEeqnarray}
and
\begin{IEEEeqnarray}{LLL}
\sigma_n^2&=&\boldsymbol{\delta}^\top\mathbf{H}_\ast^\top\boldsymbol{\Sigma}^{-1}\boldsymbol{\Sigma}_\ast\boldsymbol{\Sigma}^{-1}\mathbf{H}_\ast\boldsymbol{\delta}.
\end{IEEEeqnarray}
The minimum probability of error is given then by
\begin{IEEEeqnarray}{LLL}
P_{e}(\boldsymbol{\delta})&=&\frac{1}{2}\left(Q\left(\frac{Z(\boldsymbol{\delta})}{\sigma_n}\right)+Q\left(\frac{Z(-\boldsymbol{\delta})}{\sigma_n}\right)\right)
\label{eq:Peh}
\end{IEEEeqnarray}
The argument of the Q-function
$\frac{Z(\boldsymbol{\delta})}{{\sigma_n}}$ is given by the
following
\begin{equation}
\frac{Z(\boldsymbol{\delta})}{{\sigma_n}}=\frac{\frac{1}{2}\boldsymbol{\delta}^\top\mathbf{H}_\ast^\top\boldsymbol{\Sigma}^{-1}\mathbf{H}_\ast\boldsymbol{\delta}^\top+\boldsymbol{\delta}^\top\mathbf{H}_\ast^\top\boldsymbol{\Sigma}^{-1}(\boldsymbol{\mu}-\boldsymbol{\mu}_\ast)}{\sqrt{\boldsymbol{\delta}^\top\mathbf{H}_\ast^\top\boldsymbol{\Sigma}^{-1}\boldsymbol{\Sigma}_\ast\boldsymbol{\Sigma}^{-1}\mathbf{H}_\ast\boldsymbol{\delta}}}.
\label{eq:Q_content_functions_equal_linear}
\end{equation}
Given that the probability of error is independent of $\boldsymbol{\btheta}_o$, the bound can now be obtained upon using the expression of $P_{e}(\boldsymbol{\delta})$ from (\ref{eq:Peh}) in (\ref{eq:ZZLB_vector_ind}).

Starting with this last result, in the following subsections we further analyze
the bound for scenarios where the means and covariance matrices of the noise
also match.

\subsubsection{Equal noise means}
\label{Equal_noise_means}
Let us suppose now the assumed model correctly predicts the noise
mean $\boldsymbol{\mu}$ as the correctly specified noise mean
$\boldsymbol{\mu}_\ast$. We have that
\begin{equation}
\boldsymbol{\mu}=\boldsymbol{\mu}_\ast.
\end{equation}
Under this assumption we find that
\begin{IEEEeqnarray}{LLL}
Z(\boldsymbol{\delta})&=&\frac{1}{2}\boldsymbol{\delta}^\top\mathbf{H}_\ast^\top\boldsymbol{\Sigma}^{-1}\mathbf{H}_\ast\boldsymbol{\delta}\IEEEnonumber\\
&=&Z(\boldsymbol{-\delta}).
\end{IEEEeqnarray}
The minimum error probability can be written as
\begin{IEEEeqnarray}{lll}
P_{e}(\boldsymbol{\delta})&=&\frac{1}{2}\left(Q\left(\frac{Z(\boldsymbol{\delta})}{\sigma_n}\right)+Q\left(\frac{Z(\boldsymbol{\delta})}{\sigma_n}\right)\right)=Q\left(\frac{Z(\boldsymbol{\delta})}{\sigma_n}\right).
\label{eq:Peh2}
\end{IEEEeqnarray}
The argument of the Q-function is as
\begin{equation}
\frac{Z(\boldsymbol{\delta})}{{\sigma_n}}=\frac{\boldsymbol{\delta}^\top\mathbf{H}_\ast^\top\boldsymbol{\Sigma}^{-1}\mathbf{H}_\ast\boldsymbol{\delta}}
{2\sqrt{\boldsymbol{\delta}^\top\mathbf{H}_\ast^\top\boldsymbol{\Sigma}^{-1}\boldsymbol{\Sigma}_\ast\boldsymbol{\Sigma}^{-1}\mathbf{H}_\ast\boldsymbol{\delta}}}.
\label{eq:Z1_eqmean}
\end{equation}
As in the previous case, the bound is obtained from (\ref{eq:ZZLB_vector_ind}) upon inserting the error probability in (\ref{eq:Peh2}).

An interesting special case is that of uncorrelated, identically
distributed observations. Consider for instance the assumed and true
model being such that
\begin{eqnarray}
  \boldsymbol{\Sigma} &=& \sigma^2 \mathbf{I} \IEEEnonumber \\
  \boldsymbol{\Sigma}_\ast &=& \sigma^2_\ast \mathbf{I} ~,
\end{eqnarray}
\noindent respectively. By considering these models for the
covariance matrices, (\ref{eq:Z1_eqmean}) can be worked as
\begin{eqnarray}
\frac{Z(\boldsymbol{\delta})}{{\sigma_n}}&=&\frac{\boldsymbol{\delta}^\top\mathbf{H}_\ast^\top\mathbf{H}_\ast\boldsymbol{\delta}}
{2\sqrt{\boldsymbol{\delta}^\top\mathbf{H}_\ast^\top\mathbf{H}_\ast\boldsymbol{\delta}}}
\frac{1}{\sigma^2}\sqrt{\frac{\sigma^4}{\sigma^2_\ast}} \\
{}&=& \frac{1}{2 \sigma_\ast}
\sqrt{\boldsymbol{\delta}^\top\mathbf{H}_\ast^\top\mathbf{H}_\ast\boldsymbol{\delta}}
~,
\end{eqnarray}
\noindent which coincides with the case of no model mismatching. The
bound tells us that in the case of Gaussian linear systems with
\emph{i.i.d.} observations specifying appropriately or not the noise
variance is not affecting the theoretical performance of estimators.
Notice that this result does not hold for the case of
non-independent and/or non-identically distributed observations, in
which case covariance mispecification impacts the theoretical
performance of the estimators.

\subsubsection{Equal noise means and equal covariance matrices}
Now we analyze the bound for the scenario where the noise covariance
matrices also match. Given all the assumptions considered before,
this case is equivalent to have a total match between the assumed
and the correctly specified model. Let us write that having
\begin{equation}
\boldsymbol{\Sigma}=\boldsymbol{\Sigma}_\ast,
\end{equation}
the minimum error probability is then found as in (\ref{eq:Peh2})
\begin{IEEEeqnarray}{LLL}
P_{e}(\boldsymbol{\delta})&=Q\left(\frac{Z(\boldsymbol{\delta})}{\sigma_n}\right)\IEEEnonumber\\
\label{eq:Peh3}
\end{IEEEeqnarray}
where $\frac{Z(\boldsymbol{\delta})}{\sigma_n}$ is now equal to
\begin{IEEEeqnarray}{LLL}
\frac{Z(\boldsymbol{\delta})}{{\sigma_n}}&=&\frac{\boldsymbol{\delta}^\top\mathbf{H}_\ast^\top\boldsymbol{\Sigma}_\ast^{-1}\mathbf{H}_\ast\boldsymbol{\delta}}{2\sqrt{\boldsymbol{\delta}^\top\mathbf{H}_\ast^\top\boldsymbol{\Sigma}_\ast^{-1}\boldsymbol{\Sigma}_\ast\boldsymbol{\Sigma}_\ast^{-1}\mathbf{H}_\ast\boldsymbol{\delta}}}\IEEEnonumber\\
&=&\frac{1}{2}\sqrt{\boldsymbol{\delta}^\top\mathbf{H}_\ast^\top\boldsymbol{\Sigma}_\ast^{-1}\mathbf{H}_\ast\boldsymbol{\delta}}.
\label{eq:Q_content_match}
\end{IEEEeqnarray}
Replacing the error probability $P_{e}(\boldsymbol{\delta})$ in
(\ref{eq:ZZLB_vector_ind}) with the one from  (\ref{eq:Peh3}) yields
the ZZB for this last case study. This coincides with the classical
ZZB result for linear Gaussian systems.

\section{Special Cases: noise follows a Gaussian mixture distribution} \label{sec:special_Gaussian_mixture}

In this section we particularize the bound for the case where the
underlying noise distribution $p_{\ast}({\mathbf{n}_\ast})$ follows
a Gaussian mixture, that is, we are considering a mismatch on the noise distribution. This situation is particularly important in
cases where data is generated from two or more populations mixed in
varying proportions. Moreover, the significance of this particular
case in this work is in the fact that Gaussian mixture distributions
can approximate arbitrarily well any given density
\cite{Feller66,Sorenson71}. Therefore, this case potentially covers
all distributions $p_{\ast}({\mathbf{n}_\ast})$, when an appropriate
approximation in terms of a Gaussian mixture model is available.
Then, in this section we deal with the case of
\begin{equation}
   p_{\ast}(\mathbf{n}_\ast) = \sum_{i=1}^L \omega_i p_{i\ast}(\mathbf{n}_\ast)   =\sum_{i=1}^L \omega_i ~ \mathcal{N}(\mathbf{n}_\ast;~ \bm{\mu}_{i\ast} , \bm{\Sigma}_{i\ast}) ~,
	\label{eq:p_mixture}
\end{equation}
\noindent where $L$ is the number of mixing components and
$\omega_i$ the $i$-th mixing coefficient such that $\sum_i
\omega_i=1$ and $\omega_i \geq 0$, $\forall i$.

The goal is now to determine which is the impact of a Gaussian mixture noise 
in the derivation of the bound. Let us remind the expression of the noise variable
$n$ given in (\ref{eq:n_add}):
\begin{equation}
n=\mathbf{n}_\ast^\top\boldsymbol{\Sigma}^{-1}(\mathbf{h}(\boldsymbol{\theta}_o)-\mathbf{h}(\boldsymbol{\theta}_o+\boldsymbol{\delta})).
\end{equation}
Note that the variable $n$ also follows a Gaussian mixture. The distribution of $n$ is then as $p_n(n)=\sum_{i=1}^L \omega_i p_{n_i}(n)$, where $p_{n_i}(n)$ is associated with $p_{i\ast}(\mathbf{n}_\ast)$. The probability in (\ref{eq:Qs1}) can be computed as
\begin{IEEEeqnarray}{LLL}
\mbox{Pr}\left(S(\boldsymbol{\btheta}_o,\boldsymbol{\delta})+n<0\right)&=&\sum_{i=1}^L \omega_i Q\left(\frac{S(\boldsymbol{\btheta}_o,\boldsymbol{\delta})+\mu_{n_i}}{\sigma_{n_i}}\right)
\label{eq:Qs1_mixture}
\end{IEEEeqnarray}
where
\begin{IEEEeqnarray}{LLL}
\mu_{n_i}=\boldsymbol{\mu}_{i\ast}^\top\boldsymbol{\Sigma}^{-1}(\mathbf{h}(\boldsymbol{\btheta}_o)-\mathbf{h}(\boldsymbol{\btheta}_o+\boldsymbol{\delta})).
\label{eq:meann_mixture}
\end{IEEEeqnarray}
and
\begin{equation}
\sigma_{n_i}^2=(\mathbf{h}(\boldsymbol{\btheta}_o)-\mathbf{h}(\boldsymbol{\btheta}_o+\boldsymbol{\delta}))^\top\boldsymbol{\Sigma}^{-1}\boldsymbol{\Sigma}_{i\ast}\boldsymbol{\Sigma}^{-1}(\mathbf{h}(\boldsymbol{\btheta}_o)-\mathbf{h}(\boldsymbol{\btheta}_o+\boldsymbol{\delta})).
\end{equation}

The probability in (\ref{eq:Qs2}) can be computed analogously,
and minimum probability of error in (\ref{eq:Pe_gaussian}) can be written as
\begin{IEEEeqnarray}{l}
P_{e}(\boldsymbol{\btheta}_o,\boldsymbol{\btheta}_o+\boldsymbol{\delta})=\\
\IEEEnonumber \sum_{i=1}^L \frac{\omega_i}{2}\left(Q\left(\frac{S(\boldsymbol{\btheta}_o,\boldsymbol{\delta})+\mu_{n_i}}{\sigma_n}\right)+Q\left(\frac{-S(\boldsymbol{\btheta}_o+\boldsymbol{\delta},\boldsymbol{\delta})-\mu_{n_i}}{\sigma_n}\right)\right).
\label{eq:Pe_gaussian_mixture}
\end{IEEEeqnarray}
\section{Computer simulations}
\label{sec:simus}
In this section we provide some examples for the special cases
shown in Sections \ref{sec:special_Gaussian} and
\ref{sec:special_Gaussian_mixture}.
\subsection{Example 1}
Let us start with a scalar parameter estimation problem.
We consider a correctly specified model for $\mathbf{x}$ defined by
\begin{IEEEeqnarray}{LLL}
\mathbf{x}&=&\mathbf{h}_\ast\theta+\mathbf{n}_\ast\IEEEnonumber\\
&=&\mathbf{h}_\ast\theta+\mathbf{n}+\mathbf{n}_c
\label{eq:example_correct}
\end{IEEEeqnarray}
where $\mathbf{h}_\ast$ is the true linear function relating $\theta$ to the observations,
and $\mathbf{n}_\ast=\mathbf{n}+\mathbf{n}_c$ is the true noise process consisting of
a zero mean multivariate additive Gaussian noise term $\mathbf{n}$ plus another zero mean multivariate additive Gaussian
 noise term $\mathbf{n}_c$. As both noise processes are zero mean, the correctly
specified noise $\mathbf{n}_\ast$ is also zero mean, $\boldsymbol{\mu}_\ast=\mathbf{0}$.
The noise terms $\mathbf{n}$ and $\mathbf{n}_c$ are considered  independent, and hence,
the covariance matrix of $\mathbf{n}_\ast$ is as
\begin{equation}
\boldsymbol{\Sigma}_\ast=\boldsymbol{\Sigma}+\boldsymbol{\Sigma}_c
\end{equation}
where $\boldsymbol{\Sigma}$ and $\boldsymbol{\Sigma}_c$ are the covariance matrices
of $\mathbf{n}$ and $\mathbf{n}_c$, respectively.

Let us now consider that the assumed additive Gaussian model for the observations $\mathbf{x}$ is
\begin{equation}
\mathbf{x}=\mathbf{h}_\ast\theta+\mathbf{n}.
\label{eq:example_assumed}
\end{equation}
Note that this model is missing the noise term $\mathbf{n}_c$ with respect to the
correctly specified model. This misspecified model is referred to as $\mathcal{M}_1$ hereafter.

The next step towards the derivation of the bound is to compute the minimum probability of error $\mathbb{P}_e$ for the described
misspecified and assumed models. Given that the assumed and misspecified noise
terms follow a Gaussian distribution, the functions of $\theta$ are equal
and linear, and the noise means are also equal (zero in this case), the expression $\mathbb{P}_e$
to be considered is the one provided in
Section \ref{Equal_noise_means}. Remember that the probability of error
for scalar parameter can be obtained substituting the vector parameters
$\boldsymbol{\theta}_o$ and $\boldsymbol{\delta}$, with $\theta_o$
and $h$, respectively. We have that
\begin{IEEEeqnarray}{LLL}
\mathbb{P}_{e}(h)&=&Q\left(\frac{Z(h)}{\sigma_n}\right)
\end{IEEEeqnarray}
where
\begin{IEEEeqnarray}{LLL}
\frac{Z(h)}{\sigma_n}&=&\frac{h\mathbf{h}_\ast^\top\boldsymbol{\Sigma}^{-1}\mathbf{h}_\ast h}
{2\sqrt{h\mathbf{h}_\ast^\top\boldsymbol{\Sigma}^{-1}\boldsymbol{\Sigma}_\ast\boldsymbol{\Sigma}^{-1}\mathbf{h}_\ast h}}\IEEEnonumber\\
&=&\frac{\mathbf{h}_\ast^\top\boldsymbol{\Sigma}^{-1}\mathbf{h}_\ast h}
{2\sqrt{\mathbf{h}_\ast^\top\boldsymbol{\Sigma}^{-1}(\boldsymbol{\Sigma}+\boldsymbol{\Sigma}_c)\boldsymbol{\Sigma}^{-1}\mathbf{h}_\ast}}\IEEEnonumber\\
&=&\frac{\mathbf{h}_\ast^\top\boldsymbol{\Sigma}^{-1}\mathbf{h}_\ast h}
{2\sqrt{\mathbf{h}_\ast^\top\boldsymbol{\Sigma}^{-1}\mathbf{h}_\ast+\mathbf{h}_\ast^\top\boldsymbol{\Sigma}^{-1}\boldsymbol{\Sigma}_c\boldsymbol{\Sigma}^{-1}\mathbf{h}_\ast}}.
\end{IEEEeqnarray}
For the sake of convenience, we define the constant
\begin{equation}
\gamma= \frac{\mathbf{h}_\ast^\top\boldsymbol{\Sigma}^{-1}\mathbf{h}_\ast}
{2\sqrt{\mathbf{h}_\ast^\top\boldsymbol{\Sigma}^{-1}\mathbf{h}_\ast+\mathbf{h}_\ast^\top\boldsymbol{\Sigma}^{-1}\boldsymbol{\Sigma}_c\boldsymbol{\Sigma}^{-1}\mathbf{h}_\ast}}
\end{equation}
so that the content of the Q-function can be expressed as $\frac{Z(h)}{\sigma_n}=\gamma h$.
The bound for the scalar parameter estimation for a $\mathbb{P}_{e}$ independent of $\theta_o$ is given in
(\ref{eq:ZZLB_scalar_ind}) as
\begin{equation}
\mbox{ZZB}=\frac{1}{T}\int_{0}^{T}h(T-h)Q\left(\gamma h\right)dh.
\end{equation}
When the argument of the Q-function above depends linearly on $h$, the integral can be solved integrating by parts
\begin{IEEEeqnarray}{LLL}
\mbox{ZZB}&=&\frac{T}{6}Q\left(T \gamma \right)+\frac{1}{4\gamma^2}\Gamma_{3/2}\left(\frac{T^2\gamma^2}{2}\right)\IEEEnonumber\\
&&-\frac{2}{3T\sqrt{2\pi}\gamma^3}\Gamma_{2}\left(\frac{T^2\gamma^2}{2}\right)
\label{eq:ZZB_by_parts}
\end{IEEEeqnarray}
where $\Gamma_a(x)$ is the incomplete gamma function given by
\begin{equation}
\Gamma_a(x)=\frac{1}{\Gamma(a)}\int_0^x e^{-v}v^{a-1}dv
\end{equation}
and $\Gamma(3/2)=\sqrt{\pi}/2$. For an interval of the \textit{a priori} distribution of $\theta$, $p_\theta(\theta)$, satisfying  $T>>\frac{1}{2\gamma}$, the bound reduces to
\begin{IEEEeqnarray}{LLL}
\mbox{ZZB}&=&\frac{1}{4\gamma^2}=\frac{{\mathbf{h}_\ast^\top\boldsymbol{\Sigma}^{-1}\mathbf{h}_\ast+\mathbf{h}_\ast^\top\boldsymbol{\Sigma}^{-1}\boldsymbol{\Sigma}_c\boldsymbol{\Sigma}^{-1}\mathbf{h}_\ast}}{(\mathbf{h}_\ast^\top\boldsymbol{\Sigma}^{-1}\mathbf{h}_\ast)^2}.
\label{eq:ZZB_assumed}
\end{IEEEeqnarray}
The condition $T>>\frac{1}{2\gamma}$ implies that the RMSE range of $\hat{\theta}$ is in the support of the \textit{a priori} distribution of the parameter $\theta$. This condition is typically satisfied by choosing a large enough value for $T$. If $T$ is not large enough, then the bound given in (\ref{eq:ZZB_by_parts}) cannot be simplified, as the \textit{a priori} support of $\theta$ cuts the domain of its estimates, $\hat{\theta}$. The quasi ML estimator for the assumed model in (\ref{eq:example_assumed}) can be obtained as described in (\ref{eq:quasiMLE}).

One can also consider a misspecified model that has the knowledge on the noise term $\mathbf{n}_c$, but not on the noise term $\mathbf{n}$. The assumed additive model for the observations is then $\mathbf{x}=\mathbf{h}_\ast\theta+\mathbf{n}_c$.
This assumed model is termed as $\mathcal{M}_2$. The ZZB for this assumed model can be easily obtained, accounting for the symmetry of the problem, analogously to the previous case, in which case one obtains the expression in (\ref{eq:ZZB_assumed}) exchanging $\boldsymbol{\Sigma}$ with $\boldsymbol{\Sigma}_c$ and viceversa.

For the sake of completeness, the bound for matching models is also considered. The
assumed  model for $\mathbf{x}$ is now equal to the true model in (\ref{eq:example_correct}):
\begin{equation}
\mathbf{x}=\mathbf{h}_\ast\theta+\mathbf{n}_\ast.
\label{eq:example_assumed3}
\end{equation}
The ZZB can be obtained from the result in (\ref{eq:Q_content_match}) particularized for
scalar parameter estimation. The content of the Q-function can be written as
\begin{IEEEeqnarray}{LLL}
\frac{Z(h)}{{\sigma_n}}&=&\frac{1}{2}\sqrt{\mathbf{h}_\ast^\top\boldsymbol{\Sigma}_\ast^{-1}\mathbf{h}_\ast}h.
\end{IEEEeqnarray}
Upon integrating by parts as shown in (\ref{eq:ZZB_by_parts}) and considering $T>>\frac{1}{2\gamma}$, the bound yields
\begin{IEEEeqnarray}{LLL}
\mbox{ZZB}=\frac{1}{\mathbf{h}_\ast^\top\boldsymbol{\Sigma}_\ast^{-1}\mathbf{h}_\ast}=\frac{1}{\mathbf{h}_\ast^\top\left(\boldsymbol{\Sigma}+\boldsymbol{\Sigma}_c\right)^{-1}\mathbf{h}_\ast}.
\label{eq:ZZB_assumed3}
\end{IEEEeqnarray}
This result coincides with the CRB for the model in (\ref{eq:example_correct})\cite{Kay:1993:FSS:151045}.

\begin{figure}[!t]
\centering
\centerline{\includegraphics[width=9.5cm]{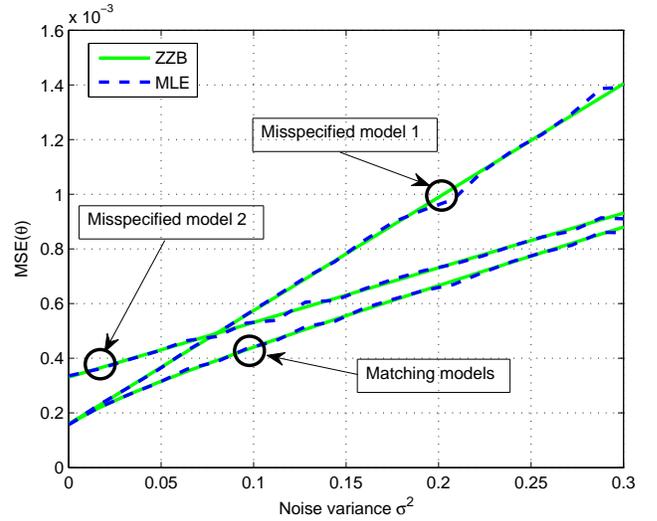}}
\caption{MSE with respect to $\sigma^2$.}
\label{example_1}
\end{figure}

Monte Carlo simulations are launched to compare the three different bounds explained above. The parameter to be estimated is set to $\theta=4$ and the length of the observation vector $\mathbf{x}$ is set to $K=500$. The noise signal $\mathbf{n}$ is additive white Gaussian noise with covariance matrix $\boldsymbol{\Sigma}=\sigma^2\mathbf{I}$. The noise variance $\sigma_n^2$ presents values ranging from $0$ to $0.3$. The noise term $\mathbf{n}_c$ is additive Gaussian noise with covariance matrix $\boldsymbol{\Sigma}_c=\sigma_c^2\mathbf{C}$, where $\sigma_c^2=0.016$ and $\mathbf{C}$ is a diagonal matrix with its diagonal elements ranging with linear spacing from $1$ to $5$. Figure \ref{example_1} illustrates the ZZBs for $\mathcal{M}_1$ derived in (\ref{eq:ZZB_assumed}) and for the analogous $\mathcal{M}_2$, as well as their corresponding MLEs. The ZZB and the MLE for the matching model in (\ref{eq:example_assumed3}) are also shown. The different bounds are tight with the corresponding estimators. Moreover, such as one would expect, $\mathcal{M}_1$ yields better estimates when the white Gaussian noise term $\mathbf{n}$ is dominant, whereas $\mathcal{M}_2$ is more robust when the variance of $\mathbf{n}$ is low, and hence, $\mathbf{n}_c$ becomes dominant. For either case, the matched model always performs equally or better. $\mathcal{M}_2$ performs equal to the matched model when $\sigma^2=0$, while $\mathcal{M}_1$ performs closer to the matched model as $\sigma^2$ increases.

\subsection{Example 2}

In this case we study the mismatch on the mean of the Gaussian noise
processes. For this purpose, we consider the following correctly specified model
for the observations $\mathbf{x}$
\begin{IEEEeqnarray}{LLL}
\mathbf{x}&=&\mathbf{h}_\ast\theta+\mathbf{n}_\ast
\label{eq:ex2_correct}
\end{IEEEeqnarray}
where $\mathbf{h}_\ast$ is the true linear function relating $\theta$ to the observations,
and $\mathbf{n}_\ast$ is a multivariate additive Gaussian noise term $\mathbf{n}$ with mean
$\boldsymbol{\mu}_\ast$ and covariance matrix $\boldsymbol{\Sigma}_\ast$. The assumed model for the observations $\mathbf{x}$ is as follows
\begin{IEEEeqnarray}{LLL}
\mathbf{x}&=&\mathbf{h}_\ast\theta+\mathbf{n}
\label{eq:ex2_assumed}
\end{IEEEeqnarray}
where $\mathbf{n}$ is additive Gaussian noise with mean  $\boldsymbol{\mu}$ and covariance
matrix $\boldsymbol{\Sigma}_\ast$. For the mismatch model under consideration, one can
find in Section \ref{functions_equal_linear} the expression of the minimum probability of error
for the bound. After particularizing the expressions in (\ref{eq:Peh}) and in (\ref{eq:Q_content_functions_equal_linear})
for scalar parameter estimation we have that
\begin{IEEEeqnarray}{LLL}
\mathbb{P}_{e}(h)&=&\frac{1}{2}\left(Q\left(\frac{Z(h)}{\sigma_n}\right)+Q\left(\frac{Z(-h)}{\sigma_n}\right)\right)
\label{eq:Peh_scalar}
\end{IEEEeqnarray}
where
\begin{IEEEeqnarray}{lll}
\frac{Z(h)}{{\sigma_n}}&=&\frac{\frac{1}{2}h\mathbf{h}_\ast^\top\boldsymbol{\Sigma}^{-1}\mathbf{h}_\ast h+h\mathbf{h}_\ast^\top\boldsymbol{\Sigma}^{-1}(\boldsymbol{\mu}-\boldsymbol{\mu}_\ast)}{\sqrt{\mathbf{h}_\ast^\top\boldsymbol{\Sigma}^{-1}\boldsymbol{\Sigma}_\ast\boldsymbol{\Sigma}^{-1}\mathbf{h}_\ast}}\IEEEnonumber\\
&=&\frac{\frac{1}{2}\mathbf{h}_\ast^\top\boldsymbol{\Sigma}^{-1}\mathbf{h}_\ast |h|+\mathbf{h}_\ast^\top\boldsymbol{\Sigma}^{-1}(\boldsymbol{\mu}-\boldsymbol{\mu}_\ast)\mbox{sgn}({h})}{\sqrt{\mathbf{h}_\ast^\top\boldsymbol{\Sigma}^{-1}\boldsymbol{\Sigma}_\ast\boldsymbol{\Sigma}^{-1}\mathbf{h}_\ast}}
\label{eq:example2_Q_content}
\end{IEEEeqnarray}
where $|\cdot|$ denotes the absolute value of its argument and $\mbox{sgn}(\cdot)$ denotes the sign of its argument.
The ZZB is obtained upon inserting the expression of $P_{e}(h)$ given in (\ref{eq:Peh_scalar}) in (\ref{eq:ZZLB_scalar_ind})
\begin{IEEEeqnarray}{LLL}
\mbox{ZZB}&=&\frac{1}{T}\int_{0}^{T}h(T-h)\frac{1}{2}Q\left(\frac{Z(h)}{{\sigma_n}}\right)dh\IEEEnonumber\\
&&+\frac{1}{T}\int_{0}^{T}h(T-h)\frac{1}{2}Q\left(\frac{Z(-h)}{{\sigma_n}}\right)dh.
\end{IEEEeqnarray}
With a change of variable, $h=-m$, the second from the expression above becomes
\begin{IEEEeqnarray}{L}
\frac{1}{T}\int_{0}^{-T}m(T+m)\frac{1}{2}Q\left(\frac{Z(m)}{{\sigma_n}}\right)dm\IEEEnonumber\\
=\frac{1}{T}\int_{-T}^{0}|m|(T-|m|)\frac{1}{2}Q\left(\frac{Z(m)}{{\sigma_n}}\right)dm.
\end{IEEEeqnarray}
The ZZB can then be expressed as follows
\begin{equation}
\mbox{ZZB}=\frac{1}{2T}\int_{-T}^{T}|h|(T-|h|)Q\left(\frac{Z(h)}{{\sigma_n}}\right)dh.
\label{finalZZB}
\end{equation}

Also, if the observations are uncorrelated and identically distributed, meaning that $\boldsymbol{\Sigma}=\sigma^2\mathbf{I}$
and $\boldsymbol{\Sigma}_\ast=\sigma^2_\ast\mathbf{I}$, the expression in (\ref{eq:example2_Q_content}) simplifies as
\begin{IEEEeqnarray}{LLL}
\frac{Z(h)}{{\sigma_n}}&=&\frac{\frac{1}{2}\mathbf{h}_\ast^\top\mathbf{h}_\ast |h|+\mathbf{h}_\ast^\top(\boldsymbol{\mu}-\boldsymbol{\mu}_\ast)\mbox{sgn}({h})}{\sqrt{\sigma_\ast^2\mathbf{h}_\ast^\top\mathbf{h}_\ast}}.
\end{IEEEeqnarray}

\begin{figure}[!t]
\centering
\centerline{\includegraphics[width=9.5cm]{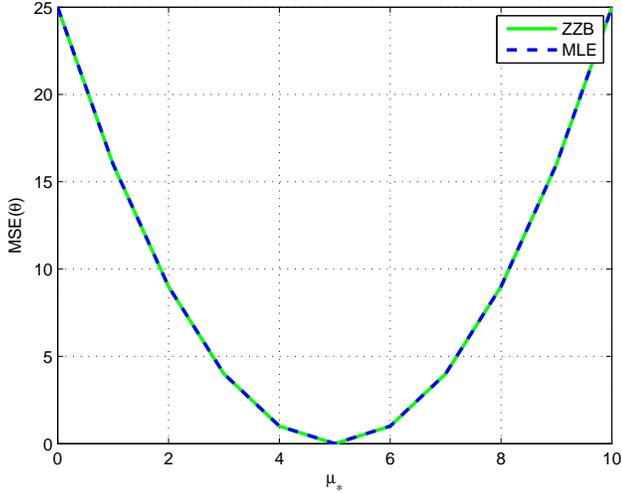}}
\caption{MSE with respect to $\mu_\ast$.}
\label{example_2}
\end{figure}

Simulations are conducted for the latter case of uncorrelated and identically distributed observations. The unknown parameter is set to $\theta=4$ and the length of the observation vector $\mathbf{x}$ is set to $K=500$. The mean of the assumed model noise process is given by $\boldsymbol{\mu}=\mathbf{1}\mu$, where $\mathbf{1}$ is the $K \times 1$ vector of ones, whereas the mean of the correctly specified model noise is $\boldsymbol{\mu}_\ast=\mathbf{1}\mu_\ast$. $\mu$ is set to $\mu=5$ and $\mu_\ast$ ranges from $0$ to $10$. The variance of the correctly specified noise is set to $\sigma^2_\ast=0.16$. Figure \ref{example_2} shows the ZZB and the MLE under this setup. The ZZB correctly predicts the bias coming from the noise mean difference.

\subsection{Example 3}

This example is devoted to the special case of Gaussian mixture
 noise analyzed in Section \ref{sec:special_Gaussian_mixture}. For this goal
we consider the following correctly specified model for the observations
$\mathbf{x}$
\begin{IEEEeqnarray}{LLL}
\mathbf{x}&=&\mathbf{h}_\ast\theta+\mathbf{n}_\ast
\end{IEEEeqnarray}
where $\mathbf{h}_\ast$ is the true linear function relating $\theta$ to the observations,
and $\mathbf{n}_\ast$ is the correctly specified noise with the following Gaussian mixture
distribution
\begin{equation}
   p_{\ast}(\mathbf{n}_\ast) = \omega_1 p_{1\ast}(\mathbf{n}_\ast)+\omega_2 p_{2\ast}(\mathbf{n}_\ast)
    \label{eq:distribution_outlier}
\end{equation}
where $\omega_1$ and $\omega_2$ are such that satisfy $\omega_1+\omega_2=1$, and  $p_{1\ast}(\mathbf{n}_\ast)$  and $p_{2\ast}(\mathbf{n}_\ast)$ are multivariate Gaussian distributions with zero mean ($\boldsymbol{\mu}_{1\ast}=\mathbf{0}$, and $\boldsymbol{\mu}_{2\ast}=\mathbf{0}$) and
covariance matrices $\boldsymbol{\Sigma}_{1\ast}$ and $\boldsymbol{\Sigma}_{2\ast}$, respectively. The distribution (\ref{eq:distribution_outlier}) is typically used to model noisy observations with outliers. In this case, the first distribution models the thermal noise and the other, with larger covariance matrix, the contribution of the outliers which occur with probability $\omega_2$. On the other side, consider the following
assumed model for the observations
\begin{IEEEeqnarray}{LLL}
\mathbf{x}&=&\mathbf{h}_\ast\theta+\mathbf{n}
\end{IEEEeqnarray}
where $\mathbf{n}$ is a zero mean Gaussian noise process with covariance matrix $\boldsymbol{\Sigma}_{1\ast}$.
One can compute the minimum probability of error $\mathbb{P}_{e}$ from the mean and the variance of the noise variable $n$ in (\ref{eq:n_add}) particularized for the model under consideration
\begin{equation}
n=-\mathbf{n}_\ast^\top\boldsymbol{\Sigma}^{-1}_{1\ast}\mathbf{h}_\ast h=\mathbf{n}_\ast^\top\mathbf{b}.
\end{equation}
where $\mathbf{b}=-\boldsymbol{\Sigma}^{-1}\mathbf{h}_\ast h$. The generic expressions for the mean and variance in mismatch Gaussian models can can be found in (\ref{eq:meanm_mixture}) and (\ref{eq:variancen_mixture_final}), respectevely. The mean is found to be
\begin{equation}
\mu_n=\omega_1 \bm{\mu}_{1\ast}^\top\mathbf{b}+\omega_2 \bm{\mu}_{2\ast}^\top\mathbf{b}=\mathbf{0}.
\end{equation}
and the variance is given by
\begin{IEEEeqnarray}{LLL}
\sigma_n^2&=&\sum_{i=1}^L \omega_i \mathbf{b}^\top\boldsymbol{\Sigma}_{i\ast}\mathbf{b}\IEEEnonumber\\
&=&\omega_1 \mathbf{h}_\ast^\top\boldsymbol{\Sigma}^{-1}_{1\ast}\mathbf{h}_\ast h^2+\omega_2 \mathbf{h}_\ast^\top\boldsymbol{\Sigma}^{-1}_{1\ast}\boldsymbol{\Sigma}_{2\ast}\boldsymbol{\Sigma}^{-1}_{1\ast}\mathbf{h}_\ast h^2.
\end{IEEEeqnarray}
The probability of error can be computed as
\begin{IEEEeqnarray}{LLL}
\mathbb{P}_{e}(h)&=&Q\left(\frac{Z(h)}{\sigma_n}\right)
\end{IEEEeqnarray}
where
\begin{equation}
\frac{Z(h)}{{\sigma_n}}=\frac{\frac{1}{2}\mathbf{h}_\ast^\top\boldsymbol{\Sigma}_{1\ast}^{-1}\mathbf{h}_\ast h}{\sqrt{\omega_1 \mathbf{h}_\ast^\top\boldsymbol{\Sigma}^{-1}_{1\ast}\mathbf{h}_\ast h^2+\omega_2 \mathbf{h}_\ast^\top\boldsymbol{\Sigma}^{-1}_{1\ast}\boldsymbol{\Sigma}_{2\ast}\boldsymbol{\Sigma}^{-1}_{1\ast}\mathbf{h}_\ast}}.
\end{equation}
As $\mathbb{P}_{e}(h)$ depends linearly on $h$, one can write  $\frac{Z(h)}{{\sigma_n}}=\gamma h$, so that the ZZB is given by $\frac{1}{4\gamma^2}$
\begin{IEEEeqnarray}{LLL}
\mathbf{ZZB}&=&\frac{\omega_1 \mathbf{h}_\ast^\top\boldsymbol{\Sigma}^{-1}_{1\ast}\mathbf{h}_\ast+\omega_2 \mathbf{h}_\ast^\top\boldsymbol{\Sigma}^{-1}_{1\ast}\boldsymbol{\Sigma}_{2\ast}\boldsymbol{\Sigma}^{-1}_{1\ast}\mathbf{h}_\ast}{(\mathbf{h}_\ast^\top\boldsymbol{\Sigma}_{1\ast}^{-1}\mathbf{h}_\ast)^2}\IEEEnonumber\\
&=&\omega_1\frac{1}{\mathbf{h}_\ast^\top\boldsymbol{\Sigma}_{1\ast}^{-1}\mathbf{h}_\ast}+(1-\omega_1)\frac{ \mathbf{h}_\ast^\top\boldsymbol{\Sigma}^{-1}_{1\ast}\boldsymbol{\Sigma}_{2\ast}\boldsymbol{\Sigma}^{-1}_{1\ast}\mathbf{h}_\ast}{(\mathbf{h}_\ast^\top\boldsymbol{\Sigma}_{1\ast}^{-1}\mathbf{h}_\ast)^2}\IEEEnonumber\\
\end{IEEEeqnarray}
for $T>>\frac{1}{2\gamma}$.

\begin{figure}[!t]
\centering
\centerline{\includegraphics[width=9.5cm]{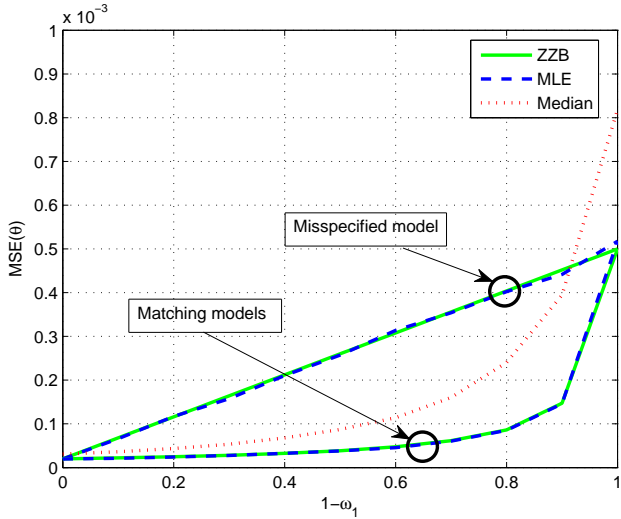}}
\caption{MSE with respect to $(1-\omega_1)$.}
\label{example_3}
\end{figure}

A test is carried out to assess the performance of the bound under a Gaussian mixture model mismatch. The unknown parameter is set as $\theta=4$ and the length of the observation vector $\mathbf{x}$ is $K=50000$. The covariance of the noise processes are given by $\boldsymbol{\Sigma}_{1\ast}=\sigma_{1\ast}^2\mathbf{I}$ and $\boldsymbol{\Sigma}_{2\ast}=\sigma^2_{2\ast}\mathbf{I}$, where $\sigma_{1\ast}=1$ and $\sigma_{2\ast}=25$. Figure \ref{example_3} illustrates the MSE as a function of $(1-\omega_1)$. The latter presents a sweeping from $0$ to $1$ in order to accommodate all possible cases. The ZZB and the MLE for the mismatch model scenario described above are shown together with their corresponding version for matching models. The performance under the misspecified model is always worse than its matching counterpart, except for the extremes cases of $\omega_1=1$ and $\omega_1=0$, where they are equivalent. Note that the example under consideration corresponds to a case study of data corrupted by outliers, that is, samples that deviate from the general distribution of the data. Outliers can easily affect the performance of classical ML estimators when the presence of such atypical observations is unknown \cite{Maronna06}. There are robust estimates that are not much influenced by outliers. One robust estimate in such models is the sample median, that is, the numerical value that separates the higher half of the observations from the lower half. For the sake of completeness, the performance of the sample median estimator is also shown. It appears from the figure that the sample median outperforms the MLE under the misspecified model, except for the limit cases.

\subsection{Example 4}

Last but not least, we assess a nonlinear multivariate parameter estimation problem. In particular, we propose a time of arrival (TOA) and amplitude estimation problem. Let us consider the following model that generates the observations $\mathbf{x}$
\begin{IEEEeqnarray}{LLL}
\mathbf{x}&=&\mathbf{h}_\ast(\btheta)+\mathbf{n}_\ast\IEEEnonumber\\
&=&\alpha\mathbf{s}_\ast(\tau)+\mathbf{n}_\ast
\end{IEEEeqnarray}
where $\mathbf{h}_\ast$ is the true function relating $\btheta$ to the observations, $\mathbf{n}_\ast$ is a zero mean multivariate additive Gaussian noise term with covariance matrix $\boldsymbol{\Sigma}_\ast$, and $\mathbf{s}_\ast$ is the correctly specified signal with energy ${E_{s_\ast}}$. The unknown parameters $\tau$ and $\alpha$ are the TOA and the received amplitude, respectively, both embedded in $\btheta=\left[\begin{array}{cc} \tau &\alpha \end{array}\right]^T$. The assumed model for the observations $\mathbf{x}$ is given by
\begin{IEEEeqnarray}{LLL}
\mathbf{x}&=&\mathbf{h}(\btheta)+\mathbf{n}_\ast\IEEEnonumber\\
&=&\alpha\mathbf{s}(\tau)+\mathbf{n}_\ast
\end{IEEEeqnarray}
here $\mathbf{h}$ is the assumed function relating $\btheta$ to the observations, and $\mathbf{s}$ is the assumed signal with energy $E_s$. Note that the inaccuracy between the models is due to different received signals.

The bound can be derived from the general case assessed in Section \ref{sec:special_Gaussian}. The variance $\sigma_n^2$ in (\ref{eq:variancen}) can be rewritten as
\begin{equation}
\sigma_n^2=(\mathbf{h}(\boldsymbol{\btheta}_o)-\mathbf{h}(\boldsymbol{\btheta}_o+\boldsymbol{\delta}))^\top\boldsymbol{\Sigma}_\ast^{-1}(\mathbf{h}(\boldsymbol{\btheta}_o)-\mathbf{h}(\boldsymbol{\btheta}_o+\boldsymbol{\delta}))
\end{equation}
as the assumed and the true noise are equal. The variable $Z(\boldsymbol{\btheta},\boldsymbol{\delta})$ in (\ref{eq:Zbthetadelta}) can be reexpressed as
\begin{IEEEeqnarray}{lll}
Z(\boldsymbol{\btheta},\boldsymbol{\delta})&=&\frac{1}{2}(\mathbf{h}(\boldsymbol{\theta}_o+\boldsymbol{\delta}))^\top\boldsymbol{\Sigma}^{-1}(\mathbf{h}(\boldsymbol{\theta}_o+\boldsymbol{\delta}))\IEEEnonumber\\
&&-\frac{1}{2}(\mathbf{h}(\boldsymbol{\theta}_o))^\top\boldsymbol{\Sigma}^{-1}(\mathbf{h}(\boldsymbol{\theta}_o))\IEEEnonumber\\
&&+(\mathbf{h}_\ast(\boldsymbol{\theta}))^\top\boldsymbol{\Sigma}^{-1}(\mathbf{h}(\boldsymbol{\theta}_o)-\mathbf{h}(\boldsymbol{\theta}_o+\boldsymbol{\delta}))
\end{IEEEeqnarray}
given that the noise term $\mathbf{n}_\ast$ is zero mean.
The minimum probability of error can be computed upon introducing the above expressions in (\ref{eq:Pe_gaussian}). The bound is obtained from the vector parameter expression given in (\ref{eq:ZZLB_vector}).

\begin{figure}[!t]
\centering
\centerline{\includegraphics[width=9.5cm]{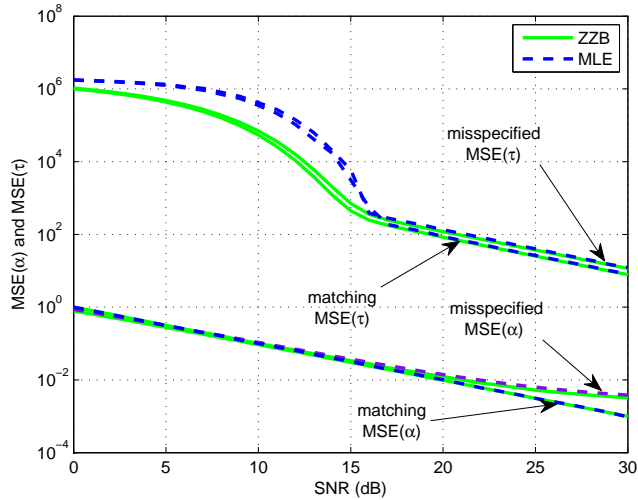}}
\caption{MSE of $\hat{\alpha}$ and $\hat{\tau}$ versusSNR.}
\label{example_4}
\end{figure}

Numerical simulations are conducted to evaluate a nonlinear parameter estimation problem in the presence of signal mismatch. The received amplitude is set to $\alpha=1$ and the length of the observation vector $\mathbf{x}$ is $K=50000$. The correctly specified signal is a triangular function with a width of $300$ samples, whereas the assumed signal is a triangular signal with a width of $200$ samples. In this case we are dealing with an inaccuracy in the signal width. Any other arbitrary discrepancy can be introduced by just choosing the appropriate waveforms. The unknown parameter $\tau$ follows discrete uniform distribution in the range $1, 2, \cdots, K$. The observations are assumed uncorrelated and identically distributed yielding a noise covariance matrix as $\boldsymbol{\Sigma}_{\ast}=\sigma_{\ast}^2\mathbf{I}$, with two-sided spectral density $\sigma_{\ast}^2=N_o/2$. Figure \ref{example_4} depicts the MSE as a function of the $\mbox{SNR}=\frac{{E_{s_\ast}}\alpha}{N_o}$. We can observe the performance of the MLE and the mismatched MLE of $\hat{\tau}$ with the corresponding bounds above in the figure. The ZZB properly anticipates the existing gap between the matching and the misspecified case. Moreover, a typical issue of non-linear estimation can be observed: the appearance of large estimation errors at low SNR. This phenomenon, known as threshold effect, takes place when, for a SNR below a certain threshold, the variance of the estimates increases considerably towards the \textit{a priori} domain of the unknown parameter \cite{Chazan1975}. The MSE performance of $\hat{\alpha}$ appears below in the figure. The quasi MLE performs close to the corresponding ZZB. In this case there is no threshold effect since amplitude estimation falls into the category of linear estimation problems.

\section{Conclusions}
\label{sec:conc}

This paper provides a methodology to derive fundamental estimation bounds in models that are misspecified, that is, cases where the assumed signal model departs from the actual phenomena being observed. Particularly, we considered the case where the assumed model is the widely used additive Gaussian model and derive the bounds for several types of model misspecification, including inaccurate functionals, noise distributions, and statistical parameters. The methodology is based on the Ziv-Zakai family of bounds. We complement the theoretical derivation of the bound with some illustrative examples, highlighting the potential of the result. Specifically, the bound is able to predict MSE of the optimal estimator, that is derived under the assumed model, when some misspecification occurs. It is noticeable that the result provides a bound on the MSE of any estimator obtained from an assumed model that departs, in some sense, from the true distribution of data. Therefore, the possible bias due to misspecification is taken into consideration, in contrast to other results based on the CRB methodology where the bias needs to be computed per estimator. This could be particularly difficult in some situations such as when no closed form expression exists for the estimator.

The bound can be potentially applied to any statistical or engineering problem where estimation relies on perfect knowledge of a signal model. Comparison of the derived bound with that derived in the usual manner without modelling errors can lead to a the sense of robustness of a method, where the furthest (in some distance sense) the MSE of the estimator lies close to the latter, the more robust it is. For instance, this was observed when assessing the performance of mean and median estimators in Figure \ref{example_3}, where the latter is known to be more robust, fact that was confirmed by our result.




\bibliographystyle{IEEEtran}

\bibliography{IEEEabrv,pau,adribiblio}


%


%




\end{document}